\PassOptionsToPackage{prologue,dvipsnames}{xcolor}
\documentclass{article}

\usepackage{arxiv}

\usepackage[utf8]{inputenc} 
\usepackage[T1]{fontenc}    
\usepackage{hyperref}       
\usepackage{url}            
\usepackage{booktabs}       
\usepackage{amsfonts}       
\usepackage{nicefrac}       
\usepackage{microtype}      
\usepackage{lipsum}		
\usepackage{graphicx}
\usepackage[square,numbers]{natbib}
\usepackage{doi}
\usepackage{tabularx}

\usepackage{amssymb,amsmath,amsthm,graphicx,tikz,caption,subcaption,pgfplots,hyperref,float}
\usepackage{siunitx, comment, multirow}
\usepackage{bm,array,arydshln}
\pgfplotsset{compat=1.18} 
\usepgfplotslibrary{groupplots}
\usepackage{pgfplotstable}
\usetikzlibrary{plotmarks,pgfplots.groupplots}
\usepackage{xstring}
\usepackage[dvipsnames]{xcolor}

\usepackage{tikz-3dplot}
\usetikzlibrary{patterns}
\usetikzlibrary{graphs, arrows.meta}
\usetikzlibrary{positioning,calc}

\usepackage{algorithm}
\usepackage{algorithmic}

\usetikzlibrary{colorbrewer}
\usepackage[normalem]{ulem}

\renewcommand{\vec}[1]{\bm{#1}}

\newcommand{\gmat}[1]{\mathbf{#1}}
\theoremstyle{definition}

\newcommand{\ASYNCWHILE}[1]{\STATE \textbf{async while} #1 \textbf{do}}
\newcommand{\ENDASYNCWHILE}{\STATE \textbf{end async while}}

\newcommand{\INDSTATE}[1][1]{\STATE\hspace{#1\algorithmicindent}}

\newcommand{\INDINDSTATE}[1][1]{\INDSTATE\hspace{#1\algorithmicindent}}

\newcommand{\INDIF}[2][1]{%
  \INDSTATE[#1]\textbf{if} #2 \textbf{then}%
}

\newcommand{\INDELSE}[1][1]{%
  \INDSTATE[#1]\textbf{else}%
}

\newcommand{\INDENDIF}[1][1]{%
  \INDSTATE[#1]\textbf{endif}%
}

\definecolor{turquoiseblue}{rgb}{0.0, 1.0, 0.94}

\pgfplotscreateplotcyclelist{exotic}{
    teal,every mark/.append style={fill=teal!80!black},mark=*\\
    orange,every mark/.append style={fill=orange!80!black},mark=square*\\
    cyan!60!black,every mark/.append style={fill=cyan!80!black},mark=triangle*\\
    red!70!white,mark=star\\
    lime!80!black,every mark/.append style={fill=lime},mark=diamond*\\
    red,densely dashed,every mark/.append style={solid,fill=red!80!black},mark=*\\
    yellow!60!black,densely dashed,
        every mark/.append style={solid,fill=yellow!80!black},mark=square*\\
    black,every mark/.append style={solid,fill=gray},mark=otimes*\\
    blue,densely dashed,mark=star,every mark/.append style=solid\\
    red,densely dashed,every mark/.append style={solid,fill=red!80!black},mark=diamond*\\
}

\usepackage[capitalise,nameinlink]{cleveref}
\Crefname{algorithm}{Algorithm}{Algorithms} 
\crefname{section}{Sec.}{Secs.}   
\Crefname{section}{Sec.}{Secs.}   
\title{A Task Parallel Orthonormalization Multigrid Method for Multiphase Elliptic Problems}

\newif\ifuniqueAffiliation

\ifuniqueAffiliation 
\author{ 
            \href{https://orcid.org/0000-0000-0000-0000}
        {\includegraphics[scale=0.06]{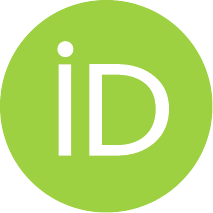} \hspace{1mm}Teoman Toprak
        }\thanks{\texttt{toprak@fdy.tu-darmstadt.de}
                \href{https://orcid.org/0000-0000-0000-0000}{Florian Kummer}
        }\\
}
\else
\usepackage{authblk}

\newbox{\orcid}\sbox{\orcid}{\includegraphics[scale=0.06]{orcid.pdf}} 
\author[1,2]{%
	\href{https://orcid.org/0009-0005-8817-5044}{\usebox{\orcid}\hspace{1mm}Teoman Toprak\thanks{\texttt{toprak@fdy.tu-darmstadt.de}}\hspace{0.5mm}},%
}
\author[1,2]{%
     \href{https://orcid.org/0000-0002-2827-7576}{\usebox{\orcid}\hspace{1mm}Florian Kummer}
}
\affil[1]{Chair of Fluid Dynamics, TU Darmstadt, Otto-Berndt-Str. 2, D-64287, Germany}
\affil[2]{Graduate School of Computational Engineering, TU Darmstadt, Dolivostr. 15, D-64293, Germany}

\fi

\renewcommand{\shorttitle}{Task Parallel Orthonormalization Multigrid Method}

\hypersetup{
pdftitle={Task Parallel Orthonormalization Multigrid Method},
pdfsubject={q-bio.NC, q-bio.QM},
pdfauthor={T. Toprak},
pdfkeywords={lienar solvers},
}

\begin{document}
\maketitle

\begin{abstract}
    Multigrid methods have been a popular approach for solving linear systems arising from the discretization of partial differential equations (PDEs) for several decades. 
    They are particularly effective for accelerating convergence rates with optimal complexity in terms of both time and space. 
    K-cycle orthonormalization multigrid is a robust variant of the multigrid method that combines the efficiency of multigrid with the robustness of Krylov-type residual minimalizations for problems with strong anisotropies. 
    However, traditional implementations of K-cycle orthonormalization multigrid often rely on bulk-synchronous parallelism, which can limit scalability on modern high-performance computing (HPC) systems.
    This paper presents a task-parallel variant of the K-cycle orthonormalization multigrid method that leverages asynchronous execution to improve scalability and performance on large-scale parallel systems.

\end{abstract}



\section{Introduction} \label{sec:Intro}

With the increasing computational power of modern hardware, expectations from scientific and engineering simulations continue to grow in both complexity and scale. 
Consequently, the numerical solution of large linear systems has become a critical component of many applications,
 especially those involving partial differential equations (PDEs), as it often represents a significant bottleneck in numerical simulations. 

Multigrid (MG) methods have long been established as a cornerstone for solving such large linear systems. 
They are particularly effective for solving linear systems arising from the discretization of elliptically dominated equations,
 such as those encountered in computational fluid dynamics, structural mechanics, and other engineering fields. 
The main advantage of multigrid methods lies in their ability to accelerate convergence rates significantly with optimal complexity in terms of both time and memory. 
Its central principle is to employ a hierarchy of discretizations, or grids, to address errors at different scales:
coarse grids efficiently reduce low-frequency error components, while fine grids handle high-frequency errors. 
A limitation of multigrid methods is that their efficiency can depend strongly on the problem structure,
 the choice of smoother, and the design of appropriate inter-grid transfer operators~\cite{Trottenberg_2001}.
For highly anisotropic, non-elliptic, or strongly heterogeneous problems, the development of effective multigrid solvers may be challenging.

In comparison, Krylov subspace methods, such as GMRES and CG, offer complementary strengths. 
They are valued for their robustness and flexibility, providing a good balance between convergence efficiency and computational cost. 
Their convergence behavior is closely linked to the spectral properties and conditioning of the system matrix, with preconditioning playing a central role in accelerating convergence~\cite{Liesen_2004}. 
In practice, Krylov methods are frequently combined with preconditioners to further enhance robustness and accelerate convergence with multigrid being one of the most prominent choices.
However, Krylov methods typically require a global communication, e.g., during inner products or matrix–vector multiplications, which becomes a bottleneck in parallel computing environments. 
In literature, there have been methods to improve the scalability of these methods, such as the pipelined GMRES \cite{Hoemmen_EECS-2010-37, Ghysels_2013, ZOU_2023_127869}, which relaxes the global communication;
 however, improvements for large-scale simulations remain limited \cite{HAWKES_2019}. 

The K-cycle multigrid method was introduced to combine the multilevel efficiency of MG with the robustness of Krylov-type residual minimization~\cite{Elman_2001,Notay_2010}.
It enhances the robustness of multigrid by
 forming several multigrid corrections and combining them through a small residual-minimization step, often implemented with GMRES or CG, into a single improved correction.
Note that, we accordingly use the term \emph{K-cycle} for MG schemes that embeds such an inner residual-minimization or Krylov-type mechanism to combine solutions obtained from different levels.
This process allows the method to better handle problems with strong anisotropies, highly variable coefficients, or other difficulties where traditional V- or W-cycles may struggle or become too expensive.
By leveraging the global error-reduction capabilities of Krylov-type residual minimization while retaining the multilevel efficiency of MG,
 the K-cycle achieves faster and more reliable convergence for challenging systems~\cite{Notay_2010}.
Alternatively, one might also use a classical Krylov method as an outer solver with multigrid as a preconditioner.

Despite these advantages, traditional K-cycle implementations still inherit the parallel scalability challenges of traditional MG and Krylov-type methods. 
Each inner minimization step requires global synchronization, introducing idle time and communication bottlenecks that become more pronounced on distributed-memory systems with high core counts. 
On coarser levels, the ratio of communication to computation increases, and the small number of degrees of freedom (DOFs) per core becomes inefficient. 
Additional global synchronization points for solution and residual updates in the Krylov-type acceleration further exacerbate these issues. 
Such factors are increasingly critical in modern high-performance computing (HPC) environments, where latency and imbalance can dominate wall-clock runtime.

Task-based parallelism offers a promising way to address these limitations. 
Similar to the additive approach, smoothing operations and coarse-grid corrections can be scheduled to overlap in time, allowing parts of the computation to proceed concurrently without fully decoupling them. 
By scheduling these tasks asynchronously, computation can be overlapped with communication, global synchronization points can be reduced, and load balance can be improved. 
Applied to the K-cycle, task-based parallelism enables concurrent execution of operations on different grid levels, coupling only the cores involved in the corresponding fine and coarse grids rather than enforcing global synchronization, while preserving the robustness of the method.
Nevertheless, this approach may introduce drawbacks, including the added complexity of task management, potential loss in convergence rates, and scheduling overhead. 
These challenges can be mitigated through careful design and implementation, allowing task parallelism to improve both scalability and runtime for large-scale simulations on modern many-core systems.

Efforts to exploit concurrency in MG have a long history, dating back at least to 1980s~\cite{McBryan_1985, GANNON_1986, Swisshelm_1986, GREENBAUM_1986, Tuminaro_1992, jones1997parallel}. 
For instance, \cite{Bastian_1998} compared an additive MG scheme enabling concurrent smoothing across grid regions, improving parallel efficiency, although multiplicative variants remain preferred for their superior convergence. 
Hybrid OpenMP--MPI strategies have also been explored, such as in \cite{Nakajima_2012,KANG_2015_2701}, where coarse-grid solves are executed on single processors to mitigate coarsest-level limitations. 

More recently, AlOnazi et al. \cite{AlOnazi_2017} proposed a hybrid MPI+OmpSs (a programming model developed at Barcelona Supercomputing Center) additive AMG that applies task-based parallelism within MPI subdomains to improve intra-node concurrency,
 but MPI ranks remains outside their task distribution, leaving inter-node operations bulk-synchronous and limiting full communication-computation overlap. 
Chaotic multigrid \cite{HAWKES_2019} eliminates global synchronization in coarser levels by allowing updates with stale data in a non-deterministic order,
 but does not provide a unified task-based framework for computation and communication. 
Moreover, the K-cycle multigrid still introduces additional synchronization points for solution and residual updates,
 which remain a scalability bottleneck at large core counts.

In this work, we present a semi-asynchronous variant of the K-cycle multigrid method, specifically the task-parallel orthonormalization multigrid method (OrthoMG). 
Our approach employs a GMRES-like residual minimization step acting on a multigrid generated search space to enhance robustness,
 while performing smoother and coarse-grid operations in parallel to improve scalability on large-scale parallel systems.
This design goes beyond existing K-cycle and synchronous multigrid schemes by leveraging task-based parallelism to overlap smoother and coarse-grid corrections across grid levels,
 reducing global synchronization and idle time.
This makes it well-suited for petascale or emerging exascale architectures,
 where minimizing synchronization and maximizing concurrency are essential for solver performance and reducing wall-clock time. 
Furthermore, its modular structure facilitates deployment on large-scale HPC systems and even enables the offloading of selected tasks to accelerators such as GPUs,
 which are increasingly prevalent in modern computing environments.

The article is organized as follows. 
\cref{sec:XDG} introduces the prototype problems, a Poisson equation with a piecewise diffusion parameter and a multiphase Stokes flow. 
Then, \cref{sec:solvers} discusses the solver and smoothers, outlining the algorithms for the respective methods. 
\cref{sec:task_parallel_orthomg} presents the task-parallel semi-asynchronous implementation of the orthonormalization K-cycle multigrid method. 
\cref{sec:Results} reports the numerical experiments with different variants and implementations along with their results. 
Finally, \cref{sec:Conclusion} provides the conclusions.


\section{Model problem definition} \label{sec:XDG}
Multigrid methods were initially applied to simple elliptic problems, but their use has since expanded to a wide range of PDEs, including parabolic and hyperbolic equations. 
In this work, we focus on the elliptic problems with heterogeneous coefficients to showcase the robustness of the K-cycle method in handling discontinuities: 
a piecewise Poisson equation and a two-phase Stokes flow.

Let us initially define the common notation and settings for both problems.
The bounded domain $\Omega \subset \mathbb{R}^d, d \in\{2,3\}$ is partitioned disjointly as:
\begin{equation}
    \Omega = \mathfrak{A} \cup \mathfrak{I} \cup \mathfrak{B},
\end{equation}
where $\mathfrak{I}$ denotes the interface separating the two bulk phases $\mathfrak{A}$ and $\mathfrak{B}$, which is assumed to be at least continuously differentiable ($C^1$).
The variables are accordingly defined piecewisely over the two phases $\mathfrak{A}$ and $\mathfrak{B}$, e.g., for a scalar $c$ and the position vector $\vec{x}$, as: 
\begin{equation}
c(\vec{x})=
    \begin{cases}
        c_{\mathfrak{A}}, & \vec{x}\in\mathfrak{A},\\[2pt]
        c_{\mathfrak{B}}, & \vec{x}\in\mathfrak{B}.
    \end{cases}
\label{eq:c_piecewise}
\end{equation}
The jump operator between bulk phases is described as $\bigl[ [c] \bigr]=(c_{ \mathfrak{B} } - c_{ \mathfrak{A} })$, pointing from $\mathfrak{A}$ to $\mathfrak{B}$. 

\subsection*{Poisson equation}
The Poisson equation with piecewise diffusion coefficient in the domain $\Omega$ is given as:
\begin{equation}
	\nabla\cdot\!\big(k \, \nabla u\big)=f \quad \text{in }\Omega,
\end{equation}
where $k$ is the predefined piecewise diffusion constant and $u$ is the scalar variable. 
The subdomains are defined by a circular interface as:
\begin{equation} \label{eq:levelset_sphere}
    \begin{aligned}
    \mathfrak{A} &:= \{ \vec{x} : |\vec{x}| < r \},\\
    \mathfrak{B} &:= \{ \vec{x} : |\vec{x}| > r \},\\
    \mathfrak{I} &:= \{ \vec{x} : |\vec{x}| = r \},
    \end{aligned}
\end{equation}
where $r$ denotes the radius of the circle in 2D and the sphere in 3D.
The interface and boundary conditions are defined as:
\begin{alignat}{4}
    \bigl[[u]\bigr] & = 0 &\quad \text{on }& \mathfrak{I}, \quad\\
    \bigl[[k\,\nabla u\cdot\vec{n}_{\mathfrak{I}}]\bigr] & = 0 &\quad \text{on }& \mathfrak{I},\\
    u & = u_\text{D} &\quad \text{on } &\partial\Omega, \quad
\end{alignat}
where $\vec{n}_{\mathfrak{I}}$ is the interface normal vector.

\subsection*{Stokes flow}
The two-phase Stokes equations, describing the motion of a viscous fluid in $\Omega$, are given as:
\begin{equation}
    \begin{split}
        \nabla\cdot (\nu \, \nabla \vec{u}) - \nabla p &= \vec{f}, \\
        \nabla \cdot \vec{u} &= 0,
    \end{split}
    \label{eq:Stokes}
\end{equation}
where $\nu$ is the piecewise constant kinematic viscosity, $\vec{u}$ is the velocity vector, $p$ is the pressure and $\vec{f}$ is the body force.
The subdomains are defined via an anisotropic (ellipsoidal) interface with a modified vector
 $\vec{x}'(\vec{x}) = ( \alpha x_0 , \, x_1)^T$ in 2D and $\vec{x}'(\vec{x}) = (\alpha x_0  , \, x_1 \, , x_2)^T$ in 3D as:
\begin{equation} \label{eq:levelset_ellipsoid}
    \begin{aligned}
    \mathfrak{A} &:= \{ \vec{x} : |\vec{x'}(\vec{x})| < r \},\\
    \mathfrak{B} &:= \{ \vec{x} : |\vec{x'}(\vec{x})| > r \},\\
    \mathfrak{I} &:= \{ \vec{x} : |\vec{x'}(\vec{x})| = r \},
    \end{aligned}
\end{equation}
where $r$ is a reference radius and $\alpha$ controls the stretching in the \(x_0\)-direction, yielding an ellipse in 2D and an ellipsoid in 3D.
The interface and boundary conditions are imposed as follows:
\begin{alignat}{4}
    \bigl[[ \vec{u} ] \bigr] & = \vec{0}       &\quad \text{on }& \mathfrak{I}, \quad\\
    \bigl[[ \big(-p\,\mathbf{I}+ \nu \left( \nabla \vec{u} + (\nabla \vec{u})^{T} \right)  \big)\vec{n}_{\mathfrak{I}} ] \bigr] & = \sigma \, \kappa \, \vec{n}_{\mathfrak{I}}        &\quad \text{on } &\mathfrak{I},\\
    \vec{u}                & = \vec{0} &\quad \text{on }& \partial\Omega, 
\end{alignat}
where $\mathbf{I}$ is the identity tensor and $\vec{n}_{\mathfrak{I}}$ is the interface normal vector. 
Moreover, $\sigma$ denotes the (constant) surface tension coefficient and $\kappa$ denotes the curvature of the interface, so that $\sigma \kappa \vec{n}_{\mathfrak{I}}$ represents the capillary surface tension force at the interface.
Note that the pressure is fixed by setting a reference value at a point in the domain to ensure uniqueness, since the homogeneous Dirichlet boundary conditions apply to the entire boundary for this problem.
\subsection*{XDG discretization}
The numerical discretization is performed using the extended discontinuous Galerkin (XDG) method, which is a special variant of the classical DG method tailored for multiphase flows. 
The XDG method employs a cut cell approach to implicitly represent the interface and allows for discontinuities in the solution across the interface. 
For the sake of brevity, we only outline the main steps of the XDG discretization, particularly the cut-cell structure. 

Following the XDG framework, the domain is discretized using a uniform background mesh $\mathfrak{K}_h$, characterized by the grid spacing $h$, consisting of regular-shaped elements $K$, which do not conform to the interface.
To account for multiple phases, background cells are described by their intersection with the subdomains $\mathfrak{s} \in \{ \mathfrak{A} , \mathfrak{B}  \}$, giving rise to \emph{phase cells} as:
\begin{equation}
    K_{i,\mathfrak{s}} := K_i \cap \mathfrak{s}, \quad \quad K_i \in \mathfrak{K}_h,
\end{equation}
whereas the corresponding extended mesh is given by $\mathfrak{K}_{h}^X := \{ K_{i,\mathfrak{s}} \}$.
Based on this construction, the XDG space is defined as:
\begin{equation}
    \mathbb{P}_p^{\text{X}}(\mathfrak{K}_h) := \{ \phi \in L^2(\Omega);
    \forall K_i \in \mathfrak{K}_h, \forall \mathfrak{s} \in \{ \mathfrak{A} , \mathfrak{B}  \} : \phi|_{K_{i,\mathfrak{s}}} \text{ is a polynomial} \nonumber \\ \text{ and deg}(\phi|_{K_{i,\mathfrak{s}}}) \leq p \}.
\end{equation}
Hence, the extended space allows for discontinuities across the interface $\mathfrak{I}$, since the basis functions are defined separately on each phase cell $K_{i,\mathfrak{s}}$.

Using the Galerkin approach, the weak formulation is obtained by multiplying the governing equation with a test function $\phi \in \mathbb{P}_p^{\text{X}}(\mathfrak{K}_h)$ and integrating over each element $K \in \mathfrak{K}_h$, performing integration by parts and introducing numerical fluxes $\hat{F}$ as:
\begin{equation}
  \sum_{K_{i,\mathfrak{s}} \in \mathfrak{K}_h^X}
  \bigl[
  \int_{K_{i,\mathfrak{s}}} (\cdots)\, \phi \, \mathrm dV  +
   \int_{K_{i,\mathfrak{s}}} (\cdots)\,\nabla \phi \, \mathrm dV  +
  \int_{\partial K_{i,\mathfrak{s}}} \hat{F}(\cdots) \,\phi \, \mathrm dS 
   \bigr] = 0
  \quad \forall \phi_m \in \mathbb{P}_p^{\text{X}}(\mathfrak{K}_h).
\end{equation}
For a detailed description and specific discretization, the reader is referred to \cite{Krause_2017, kummer_extended_2017}.

 The approximate solution is sought with the same dimension of test functions in the same polynomial space, i.e., $u_h\in \mathbb{P}_p^{\text{X}}(\mathfrak{K}_h)$, and on each phase cell $K_{i,\mathfrak{s}}$, taking the form:
\begin{equation}
    u_h|_{K_{i,\mathfrak{s}}}
  \approx \sum_{n=1}^{N} \tilde{u}_{i,\mathfrak{s},n} \,\phi_{i,\mathfrak{s},n}(\vec{x}),
\end{equation}
where $\tilde{u}_{i,\mathfrak{s},n}$ is a DG coefficient representing a degree of freedom (DOF), $N$ is the total number of local DOFs, and $\phi_{i,\mathfrak{s},n}(\vec{x})$ are the trial basis functions with $\text{supp}(\phi_{i,\mathfrak{s},n}(\vec{x})) = K_{i,\mathfrak{s}}$. 
For the Stokes problem, the velocity and pressure are approximated by $\vec{u}_h \in [\mathbb{P}_p^{\text{X}}(\mathfrak{K}_h)]^d$ and $p_h \in \mathbb{P}_q^{\text{X}}(\mathfrak{K}_h)$, respectively,
 where $q=p-1$ to satisfy the LBB condition~\cite{Babuska_1973, Brezzi_1974}.

When the phase cells are intersected by the interface $\mathfrak{I}$, i.e., if the measure of $K_i \cap \mathfrak{I}$ is positive, 
 they are referred to as cut cells, as they do not occupy the full background element volume. 
For these cells, the numerical integration method proposed by \cite{saye_high-order_2015} is employed. 
For small-cut cells, an agglomeration technique with a threshold of $0.1$ is applied to merge them with neighboring cells to ensure numerical stability \cite{Toprak_2025}.
The spatial discretization is based on the Symmetric Interior Penalty (SIP) method.
Since both problems involve only diffusive terms, no convective flux splitting is required. 

\section{Solvers and smoothers} \label{sec:solvers}
This section presents the solvers and preconditioners employed in this work. 
The orthonormalization multigrid method~\cite{Kummer_2021_BoSSS}, which forms the core of our approach, is discussed in detail. 
The Additive Schwarz and Block-Jacobi methods are then briefly described, as they serve as smoothers within the multigrid framework.
\subsection{Orthonormalization multigrid method}
The OrthoMG is a variant of the K-cycle multigrid approach that integrates Krylov-type residual minimization into the multigrid framework. 
In this variant, proposed in \cite{Kummer_2021_BoSSS}, the key idea is to augment the multigrid iteration with a residual orthonormalization
 and minimization procedure applied to both the smoother and coarse-grid correction steps, in a manner reminiscent of GMRES iteration.
\begin{algorithm}[ht!]
    \caption{Residual Minimization $(\vec{x}, \vec{r}, \gmat{W}, \gmat{Z}) := \text{RM}(\gmat{A}, \vec{r}_0, \vec{z}, \gmat{W}, \gmat{Z} )$}
    \label{alg:ResMinimi}
    \begin{algorithmic}
    \REQUIRE 
    Matrix $\gmat{A}$, initial residual $ \vec{r}_0$,
    solution $\vec{z}$,  
    an orthonormal basis  
    $\gmat{W} = (\vec{w}_1,\ldots,\vec{w}_{l}) \in \mathbb{R}^{L \times l}$  
    and its search-direction matrix $\gmat{Z}$ related to the system matrix $\gmat{A}$ by $\gmat{W} = \gmat{A} \gmat{Z}$ 
    \ENSURE An approximate solution $\vec{x}$ with residual $\vec{r}$ 
    whose $L^2$-norm does not exceed that of the initial residual from $\vec{x}_0$,
     and updated matrices $\gmat{W}, \gmat{Z} \in \mathbb{R}^{L \times (l+1)}$    
    \STATE $\vec{w}_{l+1} := \gmat{A} \vec{z}$  
    \FORALL{columns $\vec{w}_{i}$ of $\gmat{W}$}      
        \STATE $\beta := \vec{w}_i \cdot \vec{w}_{l+1}$  
        \STATE $\vec{w}_{l+1} := \vec{w}_{l+1} - \beta \vec{w}_i$  
        \STATE $\vec{z} := \vec{z} - \beta \vec{Z}_{-,i}$  
    \ENDFOR  
    \STATE $\gamma := 1 /  \|\vec{w}_{l+1}\|_2$,  
           $\vec{w}_{l+1} := \gamma \vec{w}_{l+1}$,  
           $\vec{z} := \gamma \vec{z}$  
    \STATE $\gmat{W} := ( \gmat{W}, \vec{w}_{l+1} )$,  
           $\gmat{Z} := ( \gmat{Z}, \vec{z} )$  
    \COMMENT{store vectors $\vec{w}_{l+1}$ and $\vec{z}$}  
    \STATE $\alpha_{l+1} := \vec{w}_{l+1} \cdot \vec{r}_0$,  
           $\vec{\alpha} := (\vec{\alpha}, \alpha_{l+1})$  
    \STATE $\vec{x} := \vec{x}_0 + \gmat{Z} \cdot \vec{\alpha}$  
    \COMMENT{optimized solution}  
    \STATE $\vec{r} := \vec{r}_0 - \gmat{W} \cdot \vec{\alpha}$  
    \COMMENT{residual for optimized solution}  
    \end{algorithmic}
\end{algorithm}
The orthonormalization is based on the modified Gram-Schmidt process, which allows for the construction of an orthonormal basis for the residual space. 
This basis is then used to minimize the residuals in the coarse-grid correction and smoothing steps, thereby enhancing convergence efficiency. 

The procedure, shown in \cref{alg:ResMinimi}, performs residual minimization by orthonormalizing a new search direction $\vec{w}_{l+1}$ against the existing basis $\gmat{W}$, normalizing it, and updating the associated coefficient vectors $\gmat{Z}$. 
The resulting basis is used to compute an optimized solution $\vec{x}$, whose residual norm is not greater than that of the initial guess $\vec{x}_0$. 
This process iteratively expands the search space spanned by the accumulated multigrid corrections and refines the solution.

\begin{algorithm}[ht!]
\caption{Orthonormalization Multigrid $(\vec{x}, \vec{r}) := \text{MG}(\gmat{A}, \vec{b},  \vec{x}_0) $}
\label{alg:orthomg}
\begin{algorithmic}[1]
\REQUIRE Matrix $\gmat{A}$, right-hand side $\vec{b}$ and an initial solution guess $\vec{x}_0$
\ENSURE An approximate solution $\vec{x}$ whose residual norm does not exceed that of the initial guess $\vec{x}_0$
\STATE $\gmat{Z} := \{\}$,\quad $\gmat{W} := \{\}$ \COMMENT{initialize as empty}
\STATE $\vec{r}_0 := \vec{b} - \gmat{A} \vec{x}_0,\quad \vec{r} := \vec{r}_0,\quad \vec{x} := \vec{x}_0$ \COMMENT{residual of interstitial solution}
\WHILE{$\|\vec{r}\| > \epsilon$}
    \STATE $\vec{z} := \mathrm{Smoother}(\gmat{A}, \vec{r}, \vec{0})$ \COMMENT{pre-smoother}
    \STATE $(\vec{x}, \vec{r}, \gmat{W}, \gmat{Z}) := \mathrm{RM}(\gmat{A}, \vec{r}_0, \vec{z}, \gmat{W}, \gmat{Z})$ \COMMENT{minimize residual of pre-smoother}
    \STATE $\vec{r}_\mathrm{c} := \gmat{R}^{\lambda} \vec{r}$ \COMMENT{restrict residual}
    \IF{dimension of matrix $L_c > L_{\min}$}
        \STATE $(\vec{z}_\mathrm{c}, \vec{r}_\mathrm{c}) := \mathrm{MG}(\gmat{A}^{\lambda+1},\vec{r}_\mathrm{c}, \vec{0} )$ \COMMENT{call multigrid on coarser level}
    \ELSE
        \STATE $\vec{z}_\mathrm{c} := \gmat{A}^{\lambda+1} \backslash \vec{r}_\mathrm{c}$ \COMMENT{using a sparse direct solver}
    \ENDIF
    \STATE $\vec{z} := \gmat{P}^{\lambda} \vec{z}_\mathrm{c}$ \COMMENT{prolongate coarse-grid correction}
    \STATE $(\vec{x}, \vec{r}, \gmat{W}, \gmat{Z}) := \mathrm{RM}(\gmat{A}, \vec{r}_0, \vec{z}, \gmat{W}, \gmat{Z})$  \COMMENT{minimize residual}
    \STATE $\vec{z} := \mathrm{Smoother}(\gmat{A}, \vec{r}, \vec{0})$ \COMMENT{post-smoother}
    \STATE $(\vec{x}, \vec{r}, \gmat{W}, \gmat{Z}) := \mathrm{RM}(\gmat{A}, \vec{r}_0, \vec{z}, \gmat{W}, \gmat{Z})$  \COMMENT{minimize residual of post-smoother}
\ENDWHILE
\end{algorithmic}
\end{algorithm}

\cref{alg:orthomg} outlines the orthonormalization multigrid procedure, where $\lambda$ denotes the level index.
It operates over a hierarchy of grids and consists of the standard multigrid components as smoothing, restriction, prolongation, and coarse-grid correction, augmented by the residual orthonormalization and minimization process after each update. 
It can be also applied recursively to coarser grids, yielding a fully multilevel approach in which the coarsest system is often solved with a sparse direct solver.
As the iterations are controlled by the residual-minimization procedure, the number of coarse-grid corrections per cycle is not fixed but adaptively determined by the residual reduction,
 causing the cycle to vary between V- and W-like behavior.
The algorithm iterates the solution until the residual norm is below a specified threshold $\epsilon$ or a user-defined maximum number of iterations is reached.
Note that the convergence threshold can be chosen individually for each grid level and is typically set to a larger value on coarser grids in practice.

\subsection{Additive Schwarz method}
    Schwarz methods are a class of domain decomposition techniques that were introduced by Hermann Schwarz in the 1870s~\cite{Schwarz_1870} as a mathematical tool to prove existence and uniqueness results for partial differential equations. 
    Over time, these methods evolved into practical numerical solvers for large-scale PDEs, 
    particularly in the context of parallel computing. Various formulations have been developed, 
    including overlapping and non-overlapping variants, as well as multiplicative and additive approaches. 
    In this work, we consider an additive Schwarz method due to its suitability for parallel computation, 
    and refer the reader to the review by \cite{Gander_2008Schwarz} for a detailed analysis of Schwarz methods.    
    \begin{algorithm}[h!]
        \caption{Additive Schwarz method $\vec{x}^{(k+1)} := \text{AS}(\gmat{A}, \vec{b}, \vec{x}_0) $}
        \label{alg:additive_schwarz}
        \begin{algorithmic}
            \REQUIRE Matrix $\gmat{A}$, right-hand side $\vec{b}$, 
            initial guess $\vec{x}_0$, 
            \ENSURE Approximate solution $\vec{x}$ such that $\|\gmat{A}\vec{x} - \vec{b}\| < \epsilon$, if converged
            \STATE \textbf{Initialize:} $k:=0$, $\vec{x}^{(0)} := \vec{x}_0$, $\vec{r}^{(0)} := \vec{b} - \gmat{A}\vec{x}_0$
            \WHILE{$k < k_\textrm{max} \; \And \; \|\gmat{A}x^{(k)} - \vec{b}\| \geq \epsilon$}
                \STATE $\vec{r}^{(k)} :=\vec{b} - \gmat{A}x^{(k)}$ \COMMENT{Compute global residual}
                \FOR{each subdomain $i = 1, \dots, N$ \textbf{in parallel}}
                    \STATE $\vec{r}_i := R_i \vec{r}^{(k)}$ \COMMENT{Restrict residual to subdomain}
                    \STATE Solve local problem: $\gmat{A}_i \vec{z}_i := \vec{r}_i$ \COMMENT{$\gmat{A}_i = \gmat{R}_i \gmat{A} \gmat{P}_i$ is the local subdomain matrix}
                    \STATE $\vec{e}_i := \gmat{P}_i \vec{z}_i$  \COMMENT{Prolong solution to global domain}
                \ENDFOR
                \STATE $\vec{e}^{(k)}: = \sum_{i=1}^N \vec{e}_i$ \COMMENT{Compute global correction}
                \STATE $\vec{x}^{(k+1)} := \vec{x}^{(k)} + \vec{e}^{(k)}$ \COMMENT{Update solution}
                \STATE $k := k + 1$ \COMMENT{Increment counter}
            \ENDWHILE
        \end{algorithmic}
    \end{algorithm}

   \cref{alg:additive_schwarz} outlines the additive Schwarz method for a simple block structure. 
   The domain is partitioned into subdomains with a spatial overlap of one cell using the METIS graph partitioning library~\cite{Karypis_1998METIS}, 
   which aims to balance the computational load while minimizing inter-subdomain communication. 
   Each subdomain problem is then solved locally, with the choice of local solver being flexible to accommodate problem-specific requirements (e.g., sparsity pattern, iterative or direct solvers).
   In this study, we employ PARDISO~\cite{Schenk_2002PARDISO} as the local solver due to its efficiency and robustness for sparse systems. 
   Consequently, all subdomains are processed concurrently within a loop, while only halo exchanges across cell overlaps require synchronization (in the context of smoother).
   To reduce runtime and memory costs, the LU factorization of each local system is stored and reused across successive sweeps.  
   Furthermore, the local solves are performed in single precision, while the residual computations remain in double precision, 
    which has been shown to preserve sufficient accuracy while significantly improving performance~\cite{Kummer_2021_BoSSS}.

\subsection{Block-Jacobi method}
The Block-Jacobi method is a generalization of the classical Jacobi method, where groups of unknowns are updated blockwise rather than individually~\cite{Saad_2003book}. 
The system matrix $\gmat{A}$ is partitioned into block-diagonal form, and each diagonal block $A_{ii}$ is inverted independently. 
This makes the method particularly suitable when variables within a block are strongly coupled, as often occurs in systems of PDEs with multiple degrees of freedom per cell. 
Since each block solve is independent, the method is naturally parallelizable and well-suited for modern high-performance computing architectures.

In practice, this structure has clear implications for parallelization: because block updates are computationally independent, no inter-block synchronization is required within one turn of the loop beyond forming the residual and global norms.
Hence, data access is local to each block, mapping well to distributed-memory as well as to GPU kernels, enabling high concurrency.
This high throughput often compensates for the smoother’s weaker per-iteration coupling, making Block-Jacobi attractive at scale.
\begin{algorithm}[h!]
    \caption{Block-Jacobi method $\vec{x}^{(k+1)} := \text{BJ}(\gmat{A}, \vec{b}, \vec{x}_0) $}
    \label{alg:block_jacobi}
    \begin{algorithmic}
        \REQUIRE Matrix $\gmat{A}$ with block structure $\{\gmat{A}_{ii}\}$,
         right-hand side $\vec{b}$,
          initial guess $\vec{x}_0$,
        \ENSURE Approximate solution $\vec{x}$ such that $\|\gmat{A}\vec{x} - \vec{b}\| < \epsilon$, if converged
        \STATE \textbf{Initialize:} $k:=0$, $\vec{x}^{(0)} = \vec{x}_0$
        \WHILE{$k < k_\textrm{max} \; \And \; \|\gmat{A}\vec{x}^{(k)} - \vec{b}\| \geq \epsilon$}
            \STATE $\vec{r}^{(k)} := \vec{b} - \gmat{A}\vec{x}^{(k)}$ \COMMENT{Compute residual}
            \FOR{each block $i = 1, \dots, N$ \textbf{in parallel}}
                \STATE Solve local system: $\gmat{A}_{ii} \Delta \vec{x}_i = \vec{r}_i^{(k)}$
                \STATE $\vec{x}^{(k+1)}_i := \vec{x}^{(k)}_i + \omega \Delta \vec{x}_i$ \COMMENT{Update block}
            \ENDFOR
            \STATE $k := k + 1$ \COMMENT{Increment counter}
        \ENDWHILE
    \end{algorithmic}
\end{algorithm}
\cref{alg:block_jacobi} outlines the Block-Jacobi iteration. 
Each update requires solving a small local linear system involving the block matrix $A_{ii}$, 
which can be solved exactly using direct methods or approximately using iterative techniques. 
The relaxation parameter $\omega$ can be adjusted to optimize convergence, we employ $\omega=1.0$ in this work.

The local subsystems can be solved using various techniques, depending on the size and properties of the blocks.
We employed an inverse-based approach for the local solves, which is efficient when many iterations are required and relatively small block sizes are used.
In this work, we employ LAPACK~\cite{anderson1999lapack} to compute and store the inverse of each block matrix $A_{ii}$ in advance.  
During the iteration, these precomputed inverses are applied directly to the local residuals via block sparse matrix-vector products (SpMV) multiple times per smoothing sweep (specifically five times in our experiments), 
 which eliminates repeated back substitution and reduces the computational overhead per sweep.
We also tested a factorization-based variant using a per-block PARDISO solver with cached factorizations, 
 but for our parameters and implementation, it was consistently slower than the inverse-based approach.

\section{Task-parallel OrthoMG method}
\label{sec:task_parallel_orthomg}
The key advantage of the OrthoMG method is its ability to efficiently handle large-scale problems by leveraging residual-based orthonormalization and minimization.
As the coarse-grid correction and smoothing steps are combined with the residual minimization process, the method offers greater flexibility for concurrent execution compared to traditional multigrid cycles.
An asynchronous task-parallel implementation of the OrthoMG method can be achieved by scheduling the fine-grid and coarse-grid tasks independently, 
allowing them to run concurrently on different processing units (e.g., MPI ranks or OpenMP threads).
Since these tasks are not interdependent, they can be processed in parallel, reducing idle time and improving overall scalability on many-core and heterogeneous architectures.

To loosen the coupling between the smoother and coarse-grid correction, 
the OrthoMG cycle is executed with a temporally overlapped coupling rather than a strictly sequential multiplicative scheme.
In this variant, the smoother path and the coarse path proceed concurrently on separate processor groups, each operating on a possibly stale view of the iterate/residual within the same outer iteration.
The smoother group performs the smoothing operation and residual orthonormalization/minimization, while the coarse group operates the restriction, coarse-grid solve, and prolongation. 
After both groups complete their tasks, a single synchronization step is performed to exchange the updated solution and residual information as discussed in \cref{sec:adaptiveLoops}.

\begin{algorithm}[ht!]
\caption{Semi-asynchronous orthonormalization multigrid with adaptive loops \\ $(\vec{x}, \vec{r}) := \text{AsyncMG}(\gmat{A}, \vec{b},  \vec{x}_0) $}
\label{alg:adaptive_orthomg}
\begin{algorithmic}[1]
\REQUIRE Matrix $\gmat{A}$, Right-hand side $\vec{b}$, initial guess $\vec{x}_0$
\ENSURE Approximate solution $\vec{x}$ with $\|\vec{r}\| \le \epsilon$
\STATE $\gmat{Z}\!:=\!\{\}$, $\gmat{W}\!:=\!\{\}$ \COMMENT{search spaces}
\STATE $\vec{r}_0\!:=\!\vec{b}-\gmat{A}\vec{x}_0$, $\vec{r}\!:=\vec{r}_0$
\STATE Launch two processor groups: \emph{smoother} and \emph{coarse}; set flags \texttt{done\_s}=\texttt{done\_c}=\textbf{false}
\WHILE{true}
  \ASYNCWHILE{\texttt{done\_c}=false}    \COMMENT{--- smoother group (multiplicative) ---}
  \INDSTATE $\vec{z}_\mathrm{s} := \mathrm{Smoother}(\gmat{A}, \vec{r}, \vec{0})$  \COMMENT{smoothing}
  \INDSTATE $(\vec{x},\vec{r}, \gmat{W}, \gmat{Z}) := \mathrm{RM}(\gmat{A}, \vec{r}_0,\vec{z}_\mathrm{s},\gmat{W},\gmat{Z})$ \COMMENT{residual minimization handled on smoother side}
  \INDSTATE \texttt{done\_s}\,$:=$\,\textbf{true}  
  \ENDASYNCWHILE \COMMENT{--- end of smoother group ---}
  \STATE
  \ASYNCWHILE{\texttt{done\_s}=false}   \COMMENT{--- coarse group (additive, once per cycle) ---}
  \INDSTATE $\vec{r}_\mathrm{c} := \gmat{R}^{\lambda}\vec{r}$ \COMMENT{restriction}
  \INDIF{$\dim(L_c) > L_{\min}$}
      \INDINDSTATE $(\vec{z}_\mathrm{c},\vec{r}_\mathrm{c}) := \mathrm{MG}(\gmat{A}^{\lambda+1}, \vec{r}_\mathrm{c}, \vec{0})$ \COMMENT{solve recursively}
  \INDELSE
      \INDINDSTATE $\vec{z}_\mathrm{c} := \gmat{A}^{\lambda+1} \backslash \vec{r}_\mathrm{c}$ \COMMENT{direct solver on coarsest level}
  \INDENDIF
  \INDSTATE $\vec{z}_\text{cor} := \gmat{P}^{\lambda}\vec{z}_\mathrm{c}$ \COMMENT{prolongation}
  \INDSTATE \texttt{done\_c}\,$:=$\,\textbf{true} 
  \ENDASYNCWHILE    \COMMENT{--- end of coarse group ---}
  \STATE
    \STATE \COMMENT{--- single synchronization per cycle ---}
    \IF{\texttt{done\_s} \textbf{and} \texttt{done\_c}} 
        \STATE \textbf{exchange} $\vec{z}_\text{cor}$ \COMMENT{smoother side receives coarse-grid correction}
        \STATE $(\vec{x}, \vec{r}, \gmat{W}, \gmat{Z}) := \mathrm{RM}(\gmat{A}, \vec{r}_0, \vec{z}_\mathrm{cor},\, \gmat{W},\, \gmat{Z})$ \COMMENT{apply coarse correction and update residual}
        \STATE \textbf{exchange} $\vec{r}$ \COMMENT{send updated residual to coarse group}
        \STATE reset \texttt{done\_s} = \texttt{done\_c} = \textbf{false}
        \IF{$\|\vec{r}\| \le \epsilon$}
            \STATE \textbf{break} \COMMENT{terminate}
        \ENDIF
    \ENDIF

\ENDWHILE
\end{algorithmic}
\end{algorithm}

\cref{alg:adaptive_orthomg} outlines the proposed task-parallel OrthoMG algorithm. 
It begins by initializing empty search spaces for the residual minimization process. 
Then, the two processor groups are launched and advanced independently under the aforementioned temporally overlapped, lagged scheme.
Synchronization occurs only once per cycle, minimizing communication overhead and allowing for efficient overlap of computation and communication.
The traversal follows the recursive two-grid structure, descending to coarser levels and returning upward. 
Although a single coarse-grid correction is performed per cycle, the method may revisit coarser levels in subsequent iterations since its behavior is controlled by residual.
Upon convergence, both groups terminate.
   
One of the primary advantages of this approach is the ability to assign dedicated processor groups to specific tasks, such as coarse-grid correction or smoothing.
By reducing the number of cores used for coarse-grid solves, the degrees of freedom per core increase, improving the computation-to-communication ratio and alleviating the coarse-grid bottleneck.

Task parallelism also reduces the need for global synchronization.
Synchronization is limited to smaller, specialized groups rather than enforced across all processors as discussed in \cref{sec:semi-async}.
This means that coarse-grid processors can proceed to the next coarse level while smoother processors continue their work on the current one, effectively overlapping communication with computation and hiding latency.

In addition, the approach supports adaptive task allocation at runtime, enabling adjustments based on workload variations and hardware availability (see \cref{sec:adaptiveLoops}). 
This adaptability is particularly valuable in heterogeneous environments, such as clusters with nodes of varying performance characteristics. 
Such flexibility helps with differences in processing speeds, memory hierarchies, and communication latencies between components, thereby maximizing resource utilization and scalability. 
In addition, the same principles can be applied to exploit differences between distinct component types within the system, such as leveraging GPUs or other accelerators for certain operations.

\subsection{Semi-asynchronousity}\label{sec:semi-async}
\begin{figure}[ht!]
    \centering
    \includegraphics[width=0.95\textwidth]{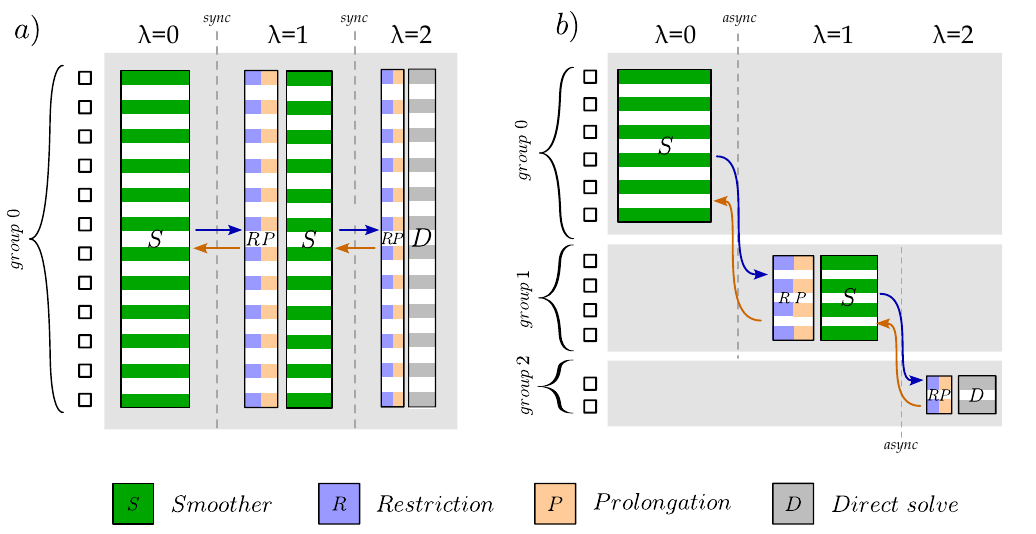}
    \caption{Schematic illustration of the synchronous (a) and task-parallel (b) orthonormalization multigrid methods with four mesh levels distributed across processing units.  
    Each column corresponds to a mesh level, while each row indicates the associated groups.  
    Symbols indicate the operations on each level, and the filled portions exhibit how work is distributed among the groups.
    \label{fig:tpSchema}}
\end{figure}
Since each residual-minimization step requires a synchronization point for updating the solution and residual, the OrthoMG method cannot be made fully asynchronous. 
Instead, it is implemented in a semi-asynchronous manner. 
In this setting, coarse-grid correction and smoothing steps can be carried out concurrently, while the solution update still enforces a synchronization point. 
Importantly, this synchronization is local to the processing units operating the involved levels, and does not extend to the entire grid hierarchy. 
Similar restrictions also appear in other asynchronous multigrid algorithms, e.g., \cite{Pou_2019}.

The synchronization step serves an essential role in our orthonormalization variant: 
the residual update itself ensures that corrections remain orthonormal, 
which by design avoids the overcorrection problem: a situation where different levels contribute nearly identical updates that would otherwise be redundant. 
In our approach, the smoother side is held responsible for residual updates, since the fine-grid data are already available there.

Both the coarse-grid correction and the smoother are executed at least once per cycle.
On the smoother side, iterations proceed multiplicatively: after each smoothing sweep,
the residual is orthonormalized as in \cref{alg:ResMinimi} and updated locally before the next sweep.
On the coarse side, the restricted right-hand side is formed once per cycle (from the state at the start of the pass), and the resulting coarse correction is applied at the end-of-cycle synchronization.
By design, it would also be possible to perform multiple coarse-grid corrections per cycle, yielding variants resembling W- or F-cycles for the OrthoMG method. 
However, such variants do not align well with either the task-parallel approach or the K-cycle, 
as they introduce additional synchronization points and undermine the efficiency of the adaptive looping between smoother and coarse-grid groups.

\subsection{Convergence criteria}\label{sec:convCriteria}
In the class of asynchronous methods, several strategies exist to control the termination and respectively convergence.
In fully asynchronous settings, imposing a global stopping test is non-trivial~\cite{HAWKES_2019}, as it requires implicit or explicit synchronization points.
A common choice is to run a fixed number of iterations per block/level and then return ~\cite{HAWKES_2019, ANZT_2013,Pou_2019}.
Another option is to dedicate a core/thread to evaluate a residual tolerance and signal termination to others.
In multigrid, the dedicated core(s) can check either only the global residual associated with the finest grid or the local residuals associated with a grid level.
Wolfson-Pou \& Chow~\cite{Pou_2019} explores both approaches on a shared-memory machine (OpenMP) and report better convergence for level-local checks, which access fresher data.
K-cycle multigrid inherently imposes residual calculation for each level when combining the solution via minimalization/orthonormalization, so level-local stopping aligns well with the algorithm.
Moreover, the synchronization is only limited between the fine and coarse grids after each iteration in the task-based parallelization.
Thus, we adopt a local-residual approach: each level checks its own residual and determines the convergence.

In our implementation, we employ this level-local termination strategy with the following explicit convergence criteria.  
The behavior of the finest grid ($\lambda = 0$) is determined by user-defined parameters with respect to the residual,
 typically in the range of $10^{-6}$ to $10^{-9}$, depending on the desired accuracy and problem characteristics.  
The first coarse grid ($\lambda = 1$) has the termination criterion based on a relative residual reduction by a factor of ten ($r \le 0.1\,r_0$) or a maximum of twenty iterations.  
The second coarse grid ($\lambda = 2$) uses a more relaxed condition, allowing at most two iterations or a residual reduction by a factor of two ($r \le 0.5\,r_0$).  
For deeper levels ($\lambda \ge 3$), an additional convergence condition is not imposed, meaning that they only perform a single iteration per visit.
The coarsest level is handled by a direct solver and is therefore excluded from the termination logic.
This configuration aligns with the observations in Notay~\cite{Notay_2010}, 
 who note that reducing the coarse-level residual by one order of magnitude or performing a small fixed number of inner Krylov iterations (3--10) 
 is sufficient for robust convergence in K-cycles.  
For consistency, we apply these criteria uniformly across all OrthoMG variants in the numerical experiments presented in the following section, 
 including both the classical and the task-parallel formulations.
\subsection{Adaptive loops}\label{sec:adaptiveLoops}
\begin{figure}[ht!]
    \centering
    \begin{tikzpicture}[font=\small,>=Latex]

\tikzset{
  box/.style={draw,rounded corners,minimum width=42mm,minimum height=10mm,align=center},
  gbox/.style={box,fill=green!25},
  pbox/.style={box,fill=violet!25},
  bbox/.style={box,fill=cyan!25},
  sbox/.style={box,fill=gray!30},
  obox/.style={box,fill=orange!30},
  ssbox/.style={draw,rounded corners,minimum width=24mm,minimum height=8mm,align=center,fill=pink!40},
  tbox/.style={draw,rounded corners,minimum width=24mm,minimum height=8mm,align=center,fill=red!40}
}
\pgfmathsetmacro{\yoff}{-2} 

\node[font=\bfseries] at (0,{6.2+\yoff}) {Smoother group};
\node[font=\bfseries] at (7,{6.2+\yoff}) {Coarse group};
\draw[->,thick] (0,{5.8+\yoff}) -- ++(0,-0.5);
\draw[->,thick] (7,{5.8+\yoff}) -- ++(0,-0.5);


\node[ssbox] at (0,{4.6+\yoff})  {Set smoother unfinished};
\draw[->,thick] (0,{4.0+\yoff}) -- ++(0,-0.5);
\node[gbox] at (0,{2.8+\yoff}) {Smoothing};
\draw[->,thick] (0,{2.2+\yoff}) -- ++(0,-0.5);
\node[pbox] at (0,{1.0+\yoff}) {Residual Minimization};



\draw[->,thick] (0,{0.4+\yoff}) -- ++(0, -0.4);
\node[ssbox] at (0,{-0.5+\yoff})  {Set smoother finished};

\draw[->,thick] (0,{-1+\yoff}) -- ++(0, -0.4);
\node[below] at (0,{-1.4+\yoff})  {Is \textbf{coarse} finished?};
\draw[->,thick] (0,{-2+\yoff}) -- ++(0, -0.4) node[right, midway] {Yes};

\draw[-,thick] (-2,{-1.6+\yoff}) -- ++(-0.5, 0) node[below, midway] {No};
\draw[-,thick] (-2.5,{-1.6+\yoff}) -- ++(0, 4.4);
\draw[->,thick] (-2.5,{2.8+\yoff}) -- ++(0.3, 0);

\node[pbox] at (0,{-3+\yoff}) {Residual Minimization};
\draw[<->] (2.5,{-3+\yoff}) -- ++(2,0)  node[midway,below=2pt]{\small Exchange};
\node[obox] at (7,{-3+\yoff}) {Prolongation};

\draw[->,thick] (0,{-3.6+\yoff}) -- ++(0, -0.5);
\node[below] at (0,{-4+\yoff}) {Did converge?};
\draw[-,thick] (-2,{-4.2+\yoff}) -- ++(-0.7, 0) node[below, midway] {No};
\draw[-,thick] (-2.7,{-4.2+\yoff}) -- ++(0, 8.7);
\draw[->,thick] (-2.7,{4.5+\yoff}) -- ++(0.5, 0);
\draw[->,thick] (0,{-4.5+\yoff}) -- ++(0, -0.5) node[right, midway] {Yes};
\node[tbox] at (0,{-5.5+\yoff}) {Terminate};

\draw[->,thick] (7,{-3.6+\yoff}) -- ++(0, -0.5);
\node[below] at (7,{-4+\yoff}) {Did converge?};

\draw[-,thick] (9,{-4.2+\yoff}) -- ++(0.7, 0) node[below, midway] {No};
\draw[-,thick] (9.7,{-4.2+\yoff}) -- ++(0, 8.7);
\draw[->,thick] (9.7,{4.5+\yoff}) -- ++(-0.5, 0);

\draw[->,thick] (7,{-4.5+\yoff}) -- ++(0, -0.5) node[right, midway] {Yes};
\node[tbox] at (7,{-5.5+\yoff}) {Terminate};


\node[ssbox] at (7,{4.6+\yoff})  {Set coarse unfinished};
\draw[->,thick] (7,{4.0+\yoff}) -- ++(0,-0.5);
\node[bbox] at (7,2.8+\yoff) {Restriction};
\draw[->,thick] (7,2.2+\yoff) -- ++(0,-0.5);
\node[sbox] at (7,1+\yoff) {Solve};


\draw[->,thick] (7,{0.4+\yoff}) -- ++(0, -0.4);
\node[ssbox] at (7,-0.5+\yoff) {Set coarse finished};

\draw[->,thick] (7,-1+\yoff) -- ++(0,-0.4);
\node[below] at (7,-1.4+\yoff) {Is \textbf{smoother} finished?};
\draw[->,thick] (7,-2+\yoff) -- ++(0,-0.4) node[right, midway] {Yes};

\draw[-,thick] (9,-1.6+\yoff) -- ++(0.5, 0) node[below, midway] {No};
\draw[-,thick] (9.5,-1.6+\yoff) -- ++(0, 2.7);
\draw[->,thick] (9.5,1.1+\yoff) -- ++(-0.3, 0);





\end{tikzpicture}
    \caption{Detailed illustration of the adaptive loop. In contrast to a standard synchronous multigrid method, the task-parallel configuration has multiple execution units, which are asynchronous but exchange information at specific points. \label{fig:adaptiveLoop}}
\end{figure}
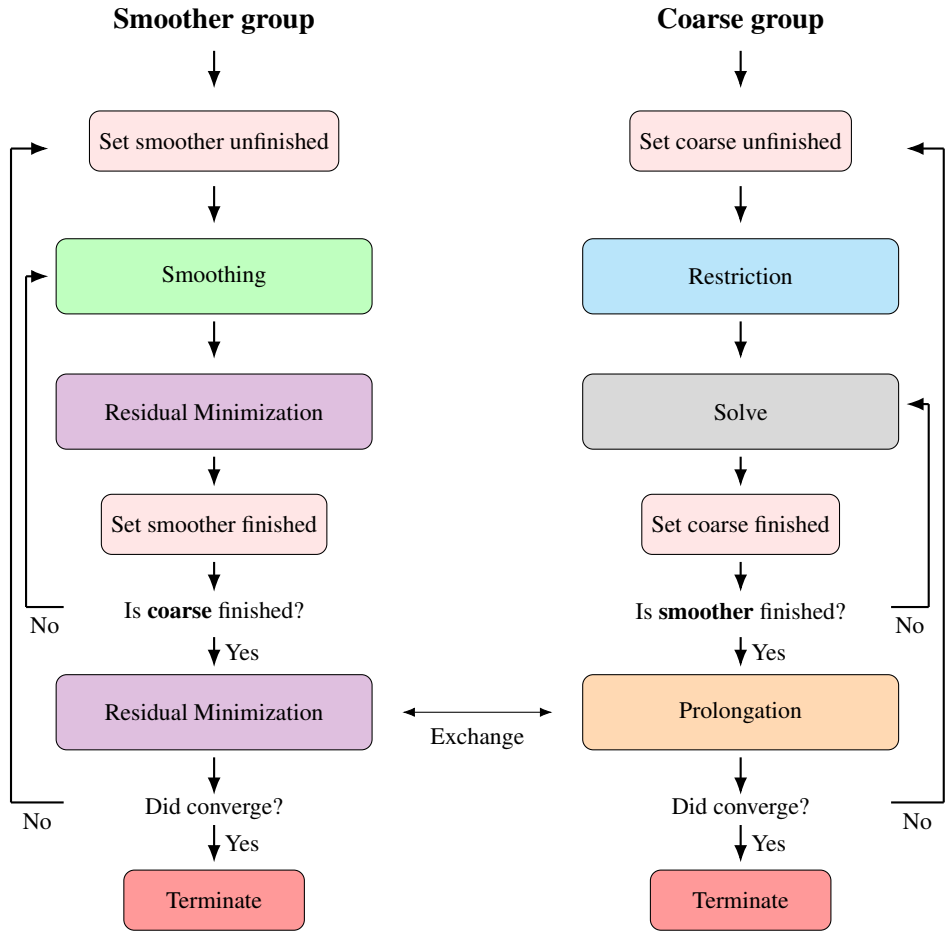
Each processor group operates in an \emph{adaptive loop}, which is \emph{theoretically chaotic} in the sense of asynchronous methods\cite{CHAZAN_1969, HAWKES_2019}, yet in practice guided by residual-based signaling rather than left fully uncoordinated.
On the smoother side, iterations proceed multiplicatively: after each smoothing step, the local residual is updated. 
On the coarse side, corrections are computed with the initial residual once per cycle.

\begin{itemize}
    \item \textbf{Signaling:} After every smoothing step (smoother group) or after producing a coarse correction (coarse group), the group sends a signal to its partner group, reporting a successful iteration.  
    \item \textbf{Waiting and continuation:} Upon sending its own signal, a group does not stop but continues sweeping with its current data until it also receives a completion signal from the other group. This ensures that neither group stalls and both remain productive even while waiting for updates.  
    \item \textbf{Exchange and restart:} Once both signals are exchanged, residuals and solutions are synchronized between groups. At this point, the coarse group may descend to the next coarser iteration (using the updated residual), while the smoother group continues on the fine grid. The adaptive loop then restarts with fresh data.  
\end{itemize}

In this way, synchronization is reduced to one exchange per cycle, after at least one smoothing step and one coarse correction. 
The loops on different mesh levels are otherwise independent, with each group working asynchronous until residual-based exchange aligns them.
Since the method is applied recursively, the coarse-grid correction has the opportunity to further reduce error components iteratively.  
If a direct solver is employed at the coarsest level, the coarse group must wait for updated residuals before proceeding, since no further iterations can be performed without new data.

\subsection{Parallelization of inter-grid operators}
A key ingredient of multigrid methods are restriction and prolongation operators, which are used to transfer residuals between different grid levels. As these operations involve both levels, there are three approaches to parallelize them:
\begin{itemize}
    \item[A] Both processor groups perform prolongation and restriction operations via sparse matrix-vector multiplication, with the results stored locally in the group that will use them. 
    \item[B] Only processors responsible for the coarse-grid level perform prolongation and restriction, receiving the restriction input and sending the prolongation output to the processors handling smoothing.
    \item[C] Only processors responsible for smoothing perform prolongation and restriction, receiving the prolongation input and sending the restriction output to the coarse-grid processors.
\end{itemize}

The first approach is the most straightforward and requires no additional matrix or data distribution operation beyond the simple domain decomposition, but it may lead to redundant synchronization and idle-time as it requires both groups. 
The second approach reduces redundancy by having only one group perform the operations, while the coarse-grid processors receive a relatively strong rectangular matrix ($a/b \approx 2^D$) for their respective levels. 
This approach can be beneficial for problems with a large number of levels, as it minimizes the amount of data transferred between groups. The third approach is presented here only for completeness; in practice, it would leave the coarse-grid processors idle, as they would not perform any prolongation or restriction operations while waiting the results.

In our numerical experiments, it is observed that the option B yields a better performance than the first option involving both levels. 
Thus, this strategy is adapted in our implementation.


\subsection{Assignment of processor groups}
In the task-parallel OrthoMG method, processor groups are assigned to either smoothing or coarse-grid correction tasks within each level.
Since the algorithm proceeds recursively through multiple grid levels, the assignment of processor groups to levels can likewise be defined recursively. 

At first glance, one might consider a coarse-to-smoother ratio depending on the spatial dimension, such as 1:3,
 to allocate more resources to the smoother tasks, which are typically more computationally intensive
and applying this ratio recursively on all levels. 
However, in practice, such ratios do not necessarily yield a good distribution of the available groups across all levels.

A more robust strategy is to base the assignment on the number of degrees of freedom (DOFs), with special treatment of the coarsest level, 
which is often solved directly and is much smaller than the finer levels. 
In this work, we follow such a heuristic: we first allocate a user-defined number of processor groups to the coarsest grid 
(in particular, a single MPI rank) to ensure that its solves remain efficient (e.g. about $10^4$ DOFs). 
Then, the remaining groups are distributed among the smoothing tasks in proportion to the DOFs on each level, 
giving priority to finer levels to reach an approximately uniform number of DOFs per core on all non-coarse levels.

\section{Numerical results} \label{sec:Results}
This section presents numerical results for the proposed task-parallel orthonormalization multigrid method applied to the model problems with a jump in diffusion coefficients
 under varying parameters. 
The method is implemented in the open-source XDG software package BoSSS\footnote{\url{https://doi.org/10.5281/zenodo.8386633}}.
The performance is evaluated in terms of convergence behavior, runtime, and parallel scalability. 
The results are compared against each other to demonstrate the advantages and disadvantages.

Four solver configurations of the orthonormalization multigrid method are considered and summarized in Table~\ref{tab:orthomg-configs}:
 \emph{additive synchronous}, \emph{multiplicative synchronous}, \emph{additive task-parallel}, and \emph{multiplicative–task-parallel hybrid}. 
 The additive synchronous variant, implemented in its standard synchronous form, only serves as the baseline for performance comparison with the task-parallel approach. 
 The multiplicative synchronous variant corresponds to the traditional multiplicative multigrid cycle with pre- and post-smoother on each level \cite{Kummer_2021_BoSSS}.
 The additive task-parallel configuration, detailed in \cref{sec:task_parallel_orthomg}, is the primary focus of this work and employs semi-asynchronous execution with an additive residual update. 
 The hybrid configuration starts with a classical synchronous multiplicative cycle on the finest level and then switches to the semi-asynchronous task-parallel cycle on the remaining levels,
  aiming to combine the robustness of the classical scheme with the improved parallel efficiency of the task-parallel variant.
  Similar to additive synchronous variant, it also serves as a baseline for performance comparison with the task-parallel approach.
 Note that we use the terms additive and multiplicative to indicate the order of residual updates,
  within the residual orthonormalization-minimalization framework the corrections are neither added nor multiplied but combined via \cref{alg:ResMinimi}.
\begin{table}[t]
  \centering
  \begin{tabularx}{\linewidth}{
    >{\arraybackslash}p{4.2cm}
    >{\arraybackslash}p{1.5cm}
    >{\arraybackslash}p{3.6cm}
    >{\arraybackslash}p{3.6cm}
    >{\arraybackslash}p{1.8cm}
  }
    \toprule
    \textbf{Configuration} & \textbf{Alg.(s)} &
    \textbf{Synchronicity} & \textbf{Residual update order} &
    \textbf{Smoothers} \\
    \midrule
    Additive synchronous &
    Alg.~\ref{alg:orthomg} &
    Synchronous &
    Additive &
    post \\

    Multiplicative synchronous &
    Alg.~\ref{alg:orthomg} &
    Synchronous &
    Multiplicative &
    pre + post \\

    Additive task-parallel &
    Alg.~\ref{alg:adaptive_orthomg} &
    Semi-asynchronous &
    Additive &
    post \\

    Hybrid\newline(multiplicative-task-parallel) &
    Algs.~\ref{alg:orthomg}, \ref{alg:adaptive_orthomg} &
    Synchronous on finest;\newline semi-asynchronous on\newline others &
    Multiplicative on finest;\newline additive on others &
    pre (finest)\newline + post (all) \\
    \bottomrule
  \end{tabularx}
  \caption{Overview of the considered orthonormalization multigrid configurations.
   The synchronicity indicates whether operations are executed strictly sequentially or can overlap across tasks.
   The residual update order distinguishes additive and multiplicative schemes: in the additive case smoother and coarse corrections are computed from the same initial residual and combined at the end,
   whereas in the multiplicative case the residual is updated after each correction.
   The smoother column indicates whether pre- or post-smoothing is employed.}
  \label{tab:orthomg-configs}
\end{table}


All numerical experiments were carried out on the Lichtenberg high-performance computing system of the TU Darmstadt under fixed hardware constraints to ensure reproducibility. 
The CPU frequency was pinned to 2.3 GHz in performance mode. 
Computations were executed exclusively on dedicated nodes supporting the AVX-512 instruction set, ensuring consistent vectorization capabilities and preventing interference from other users.
In addition, the memory allocation was fixed to 3.8 GB per CPU core. 
Following the memory-related observations for BoSSS~\cite{Kummer_2021_BoSSS}, 
 we restrict the maximum number of degrees of freedom per core to approximately $25 \times 10^3$ in all experiments.

Distributed memory parallelism is achieved via the Message Passing Interface (MPI) 
 and the task parallel implementation described in \cref{sec:task_parallel_orthomg} is realized through MPI communicator groups with non-blocking communication between them.
For each MPI rank, six threads are employed, 
  ensuring that the number of requested cores and threads are matched.
Within each MPI rank, shared memory parallelism is achieved via either OpenMP in third-party libraries (e.g., BLAS, LAPACK, PARDISO) or Task Parallel Library (TPL) feature of \emph{.NET}, on which BoSSS is implemented.
In particular, block solves in smoothers are parallelized with TPL.
MPI ranks were bound to physical cores with NUMA-aware allocation and hyper-threading was disabled.

All solvers are executed under identical conditions, including problem size, discretization parameters and domain decomposition, in order to ensure comparability.
The convergence criterion is defined as a relative residual reduction of $10^{-8}$ with respect to the initial residual norm, or reaching an absolute tolerance of $10^{-8}$.
The reported runtimes correspond to wall-clock times of the solvers, averaged over multiple runs (3-4 times) to reduce the impact of fluctuations.
Scaling behavior is analyzed by varying both the number of cores and the problem size.

\subsection{Benchmark problems} \label{sec:benchmarks}
We consider two representative benchmark problems with respect to \cref{sec:XDG}, described on a $d$-dimensional hypercubic domain $[-L,L]^d$.

For the Poisson problem with discontinuous diffusion coefficients, the interface is a circular with the radius $r = 0.7L$ centered in the domain as defined in \cref{eq:levelset_sphere}. 
A jump in the diffusion coefficient is imposed from inside to outside the shape, with $k_\mathfrak{B}/k_\mathfrak{A} = 1000$, and a constant right-hand side $f = 1$ is prescribed. 
For all boundary faces, homogeneous Dirichlet conditions are applied.

Similarly, the Stokes flow problem with discontinuous viscosity is tested.
The interface is defined to create an ellipsoidal shape of radius $r=0.49L$, centered in the domain as described by \cref{eq:levelset_ellipsoid}
 with a stretching factor $\alpha=2$.
A jump in the viscosity is imposed from outside to inside the ellipse with respect to the physical parameters of air/water pair leading to 
$\nu_\mathfrak{A}/\nu_\mathfrak{B}=14.25$, and the surface tension calculated with respect to $\sigma=0.072 \mathrm{ N/m}$.
No additional body forces are considered and homogeneous Dirichlet conditions are applied on all boundary faces.
Additionally, a minimum number of smoothing steps is set to five for the Stokes problem (applied to all solver configurations) to ensure sufficient smoothing and compatibility,
 as the task-parallel approach may already lead to multiple steps due to the adaptive loops described in \cref{sec:adaptiveLoops}.

The main parameters of both benchmark problems are summarized in 
Table~\ref{tab:benchmarks}.

\begin{table}[htbp]
  \centering
  \begin{tabular}{lcccc}
    \toprule
    Problem & Interface & Jump & Boundary conditions & Additional parameters \\
    \midrule
    Poisson & circle, $r = 0.7L$ &
      $k_\mathfrak{B}/k_\mathfrak{A} = 1000$ &
      homogeneous Dirichlet &
      $f = 1$ \\
    Stokes  & ellipse, $r = 0.49L$, $\alpha = 2$ &
      $\nu_\mathfrak{A}/\nu_\mathfrak{B} = 14.25$ &
      homogeneous Dirichlet  &
      $\sigma = 0.072$ [N/m], $\vec{f} = \vec{0}$ \\
    \bottomrule
  \end{tabular}
  \caption{Summary of benchmark problems and main parameters.}
  \label{tab:benchmarks}
\end{table}

\subsection{Choice of smoother} \label{sec:smootherComp}
\begin{figure}[h!]
    \centering
    \begin{center}
\begin{tikzpicture}
\begin{groupplot}[
  group style={
    group size=2 by 2,
    horizontal sep=18mm, vertical sep=16mm,
    xlabels at=edge bottom,
    ylabels at=edge left,
  },
  xlabel={Cores},
  ylabel={Wall-clock time [s]},
  ybar,
  enlarge x limits=0.2,
  width=6.5cm, height=5cm,
  ymode=log,
  log origin=infty, 
  ymin=1, ymax=150,
  legend style={at={(0.5,-0.15)}, anchor=north, legend columns=-1}
]

\nextgroupplot[title={(a) Multiplicative synchronous}, symbolic x coords={48, 192, 768}, xtick=data, legend to name=SharedLegend]

\addplot+[blue, fill=blue!50] coordinates {(48,91.0579271) (192, 51.4325875) (768, 96.5450451)};
\addlegendentry{BlockJacobi}

\addplot+[red, fill=red!50] coordinates {(48, 7.453119) (192, 2.5104713) (768, 5.1139766)};
\addlegendentry{Schwarz}

%
%

\nextgroupplot[title={(b) Additive task-parallel}, symbolic x coords={48, 192, 768}, xtick=data]

\addplot+[blue, fill=blue!20] coordinates {(48, 62.1616504) (192, 28.273288) (768, 6.7514704)};

\addplot+[red, fill=red!20] coordinates {(48, 2.158159) (192, 2.0815832) (768, 1.1925235)};








%
%

\end{groupplot}
\end{tikzpicture}
\end{center}

\begin{center}
\pgfplotslegendfromname{SharedLegend}
\end{center}
    \caption{Runtime comparison of OrthoMG with Block-Jacobi (blue) and Additive Schwarz (red) smoothers for varying numbers of cores for the Poisson problem with $128\times128$ cells and $p=5$. 
    The left (multiplicative synchronous) and right (additive task-parallel) panels correspond to different parallelization strategies given in \cref{tab:orthomg-configs}.
 \label{fig:smootherComp}}
\end{figure}
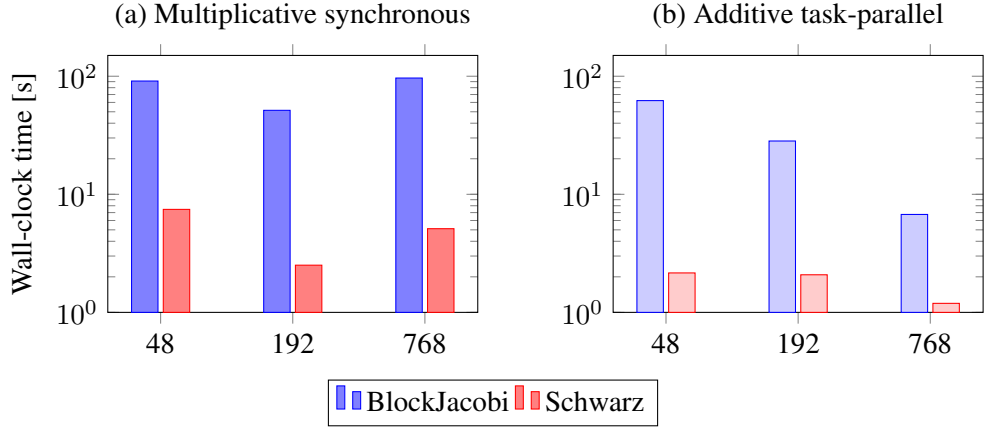
The choice of smoother plays a significant role in the overall performance and convergence behavior of multigrid methods.
This section compares the additive Schwarz (Alg.~\ref{alg:additive_schwarz}) and the Block-Jacobi (Alg.~\ref{alg:block_jacobi}) methods employed
 as smoothers within the OrthoMG framework.
Both smoothers are tested under the classical multiplicative synchronous and the proposed additive task-parallel variants as described in \cref{tab:orthomg-configs}.

We consider both benchmark problems given in \cref{sec:benchmarks} for 2D. 
The computational grids are resolved with $128\times128$ and $256\times256$ cells, while the polynomial degree is varied as $p=2, 3, 5$.
In addition, the number of cores is slightly varied to verify its influence on the performance; however, the main focus is on the comparison between the two smoothers.
The scalability behavior is analyzed by varying both the number of cores and the problem size in \cref{sec:runtimeComp}.

\cref{fig:smootherComp} illustrates the runtime behavior of the selected test cases for the two smoother configurations for the Poisson problem.
As presented, Block-Jacobi consistently leads to significantly higher runtimes than additive Schwarz for both parallelization strategy. 
We observe the same trend for all tested configurations, independent of the number of cores, grid resolution, or polynomial degree. 
Therefore, only one representative case is shown here for clarity.

This difference between smoothers can be attributed to the lack of overlap in Block-Jacobi, 
 which reduces smoothing efficiency and hence requires more iterations for convergence. 
In contrast, additive Schwarz utilizes overlapping subdomains, yielding faster convergence.
Tables~\ref{tab:poisson2d_schwarz_blockjacobi_384} summarizes the iteration counts with varying polynomial degrees and grid resolutions for the 2D Poisson problem.
For all cases, the additive Schwarz smoother demonstrates significantly lower iteration counts compared to the Block-Jacobi method,
 in accordance with the runtime results presented in \cref{fig:smootherComp}.
This confirms the superior convergence properties of the additive Schwarz method over Block-Jacobi in this context.
Moreover, the iteration counts for the additive Schwarz smoother remain relatively stable across different configurations, 
 while the Block-Jacobi method shows more variability.
\begin{table}[ht!]
    \centering
    \begin{tabular}{c c c c c c}
        \toprule
        \textbf{Grid Resolution} & \textbf{$p$} &
        \multicolumn{2}{c}{\textbf{Schwarz}} &
        \multicolumn{2}{c}{\textbf{Block Jacobi}} \\
        \cmidrule(lr){3-4}\cmidrule(lr){5-6}
        & & Multiplicative Sync & Task-parallel & Multiplicative Sync & Task-parallel \\
        \midrule
        $128\times128$ & 2 & 5   & 5   & 15   & 9   \\
        $128\times128$ & 3 & 5   & 4   & 19   & 7.5 \\
        $128\times128$ & 5 & 5   & 4   & 25   & 8.5 \\
        $256\times256$ & 2 & 6   & 5   & 15   & 9   \\
        $256\times256$ & 3 & 6   & 5   & 17   & 8.5 \\
        $256\times256$ & 5 & 6   & 5   & 38 & 21  \\
        \bottomrule
    \end{tabular}
    \caption{Average iteration counts for the 2D Poisson problem with OrthoMG, comparing the multiplicative synchronous and the additive task-parallel configurations.
    Note that due to the multithreading within each smoother or adaptive loops within the task-parallel approach, the iteration counts may vary slightly between runs.}
    \label{tab:poisson2d_schwarz_blockjacobi_384}
\end{table}

For the Stokes problem, a much more pronounced difference between the smoothers is observed.
Unlike the Poisson problem, where both smoothers yield acceptable convergence, the Block-Jacobi smoother performs substantially worse for the Stokes system.
In some configurations, it leads to almost stagnant convergence, so that only additive Schwarz remains practically viable.
This is consistent with the saddle-point structure of the Stokes system (cf.~\cref{eq:Stokes}).
Block-Jacobi damps high-frequency velocity modes but fails to reduce pressure modes and divergence errors.
On the other hand, the overlapping additive Schwarz smoother provides sufficient coupling for robust multigrid convergence.
Since we employ a monolithic formulation for the Stokes system, only the additive Schwarz smoother is hence considered in the following sections.
Table~\ref{tab:stokes2d_schwarz_blockjacobi_384} summarizes the iteration counts for the 2D Stokes problem with the additive Schwarz smoother
 under varying polynomial degrees and grid resolutions.
\begin{table}[ht!]
    \centering
    \begin{tabular}{c c c c c}
        \toprule
        \textbf{Grid Resolution} & \textbf{$p$} &
        \multicolumn{2}{c}{\textbf{Iterations}} \\
        \cmidrule(lr){3-4}
        & & Multiplicative Sync & Task-parallel \\
        \midrule
        $128\times128$ & 2 & 15 & 14 \\
        $128\times128$ & 3 & 12 &  9 \\
        $256\times256$ & 2 & 18 & 10 \\
        $256\times256$ & 3 & 19 &  10 \\
        \bottomrule
    \end{tabular}
    \caption{Average iteration counts for the 2D Stokes problem with OrthoMG, comparing the multiplicative synchronous and the additive task-parallel configurations.}
    \label{tab:stokes2d_schwarz_blockjacobi_384}
\end{table}

\subsection{Convergence behavior} \label{sec:convComp}
\begin{figure}[ht!]
    \centering
    \begin{tikzpicture}
\begin{groupplot}[
  group style={
    group size=2 by 1,
    horizontal sep=14mm, vertical sep=16mm,
    xlabels at=edge bottom,
    ylabels at=edge left,
  },
  ymode=log,
  xlabel={Iteration},
  ylabel={Residual},
  grid=both,
  legend pos=north east,
  ymin=1e-10, ymax=1e+0,      
  scale only axis, width=6cm, height=4.5cm,
]


\pgfkeys{
  /mods/.is family,
  /mods/initialMod/.initial=1,
  /mods/coarseMod/.initial=1,
  /mods/smootherMod/.initial=1
}

\pgfkeys{
  /mods/initialMod=1,
  /mods/coarseMod=1,
  /mods/smootherMod=1
}

\newcount\initialcount
\newcount\coarsecount
\newcount\smoothercount

\initialcount=0
\coarsecount=0
\smoothercount=0

\nextgroupplot[
    title={(a) Multiplicative synchronous},
    xtick={3,6,9, 12, 16},
    xticklabels={1, 2, 3, 4, 5},
]

\initialcount=0
\coarsecount=0
\smoothercount=0

\addplot+[
  only marks,
  mark=triangle*,
  mark size=2pt,
  mark options={fill=blue!70!blue,draw=black,rotate=-90},
  x filter/.code={%
    \IfStrEq{\thisrow{type}}{initial}{
        \advance\initialcount by 1
        \ifnum\initialcount=\pgfkeysvalueof{/mods/initialMod}
            \initialcount=0
        \else
            \def\pgfmathresult{}
        \fi
    }{
        \def\pgfmathresult{}
    }
  },
] table[x=step,y=residual,col sep=comma]{figures/poissonCore48_Res128_p3_Mutliplicative.txt};
\addlegendentry{Initial}

\addplot+[
  only marks,
  mark=*,
  mark size=2pt,
  mark options={fill=green!70!black,draw=black},
  x filter/.code={%
    \IfStrEq{\thisrow{type}}{smoother}{
        \advance\smoothercount by 1
        \ifnum\smoothercount=\pgfkeysvalueof{/mods/smootherMod}
            \smoothercount=0
        \else
            \def\pgfmathresult{}
        \fi
    }{
        \def\pgfmathresult{}
    }
  },
] table[x=step,y=residual,col sep=comma]{figures/poissonCore48_Res128_p3_Mutliplicative.txt};
\addlegendentry{Smoother}

\addplot+[
  only marks,
  mark=square*,
  mark size=2pt,
  mark options={fill=red!70!black,draw=black},
  x filter/.code={%
    \IfStrEq{\thisrow{type}}{coarse}{
        \advance\coarsecount by 1
        \ifnum\coarsecount=\pgfkeysvalueof{/mods/coarseMod}
            \coarsecount=0
        \else
            \def\pgfmathresult{}
        \fi
    }{
        \def\pgfmathresult{}
    }
  },
] table[x=step,y=residual,col sep=comma]{figures/poissonCore48_Res128_p3_Mutliplicative.txt};
\addlegendentry{Coarse}

\addplot+[
  only marks,
  mark=x,
  mark size=2.5pt,
  mark options={draw=black!70!black, line width=1.5pt},
  x filter/.code={%
    \IfStrEq{\thisrow{type}}{final}{}{%
      \def\pgfmathresult{}
    }%
  },
] table[x=step,y=residual,col sep=comma]{figures/poissonCore48_Res128_p3_Mutliplicative.txt};
\addlegendentry{Final}


\nextgroupplot[title={(b) Additive task-parallel},
    xtick={6,16,27,36,47},
    xticklabels={1, 2, 3, 4, 5},
    ]

\initialcount=0
\coarsecount=0
\smoothercount=0

\addplot+[
  only marks,
  mark=triangle*,
  mark size=2pt,
  mark options={fill=blue!70!blue,draw=black,rotate=-90},
  x filter/.code={%
    \IfStrEq{\thisrow{type}}{initial}{
        \advance\initialcount by 1
        \ifnum\initialcount=\pgfkeysvalueof{/mods/initialMod}
            \initialcount=0
        \else
            \def\pgfmathresult{}
        \fi
    }{
        \def\pgfmathresult{}
    }
  },
] table[x=step,y=residual,col sep=comma]{figures/poissonCore48_Res128_p3_TaskParallel.txt};
\addlegendentry{Initial}

\addplot+[
  only marks,
  mark=*,
  mark size=2pt,
  mark options={fill=green!70!black,draw=black},
  x filter/.code={%
    \IfStrEq{\thisrow{type}}{smoother}{
        \advance\smoothercount by 1
        \ifnum\smoothercount=\pgfkeysvalueof{/mods/smootherMod}
            \smoothercount=0
        \else
            \def\pgfmathresult{}
        \fi
    }{
        \def\pgfmathresult{}
    }
  },
] table[x=step,y=residual,col sep=comma]{figures/poissonCore48_Res128_p3_TaskParallel.txt};
\addlegendentry{Smoother}

\addplot+[
  only marks,
  mark=square*,
  mark size=2pt,
  mark options={fill=red!70!black,draw=black},
  x filter/.code={%
    \IfStrEq{\thisrow{type}}{coarse}{
        \advance\coarsecount by 1
        \ifnum\coarsecount=\pgfkeysvalueof{/mods/coarseMod}
            \coarsecount=0
        \else
            \def\pgfmathresult{}
        \fi
    }{
        \def\pgfmathresult{}
    }
  },
] table[x=step,y=residual,col sep=comma]{figures/poissonCore48_Res128_p3_TaskParallel.txt};
\addlegendentry{Coarse}

\addplot+[
  only marks,
  mark=x,
  mark size=2.5pt,
  mark options={draw=black!70!black, line width=1.5pt},
  x filter/.code={%
    \IfStrEq{\thisrow{type}}{final}{}{%
      \def\pgfmathresult{}
    }%
  },
] table[x=step,y=residual,col sep=comma]{figures/poissonCore48_Res128_p3_TaskParallel.txt};
\addlegendentry{Final}

\end{groupplot}
\end{tikzpicture}
    \caption{Convergence comparison of configurations of OrthoMG for the 2D Poisson problem with the additive Schwarz smoother, polynomial degree $p=3$, resolution of $128^2$ cells simulated with 48 cores \label{fig:convComp}}
\end{figure}
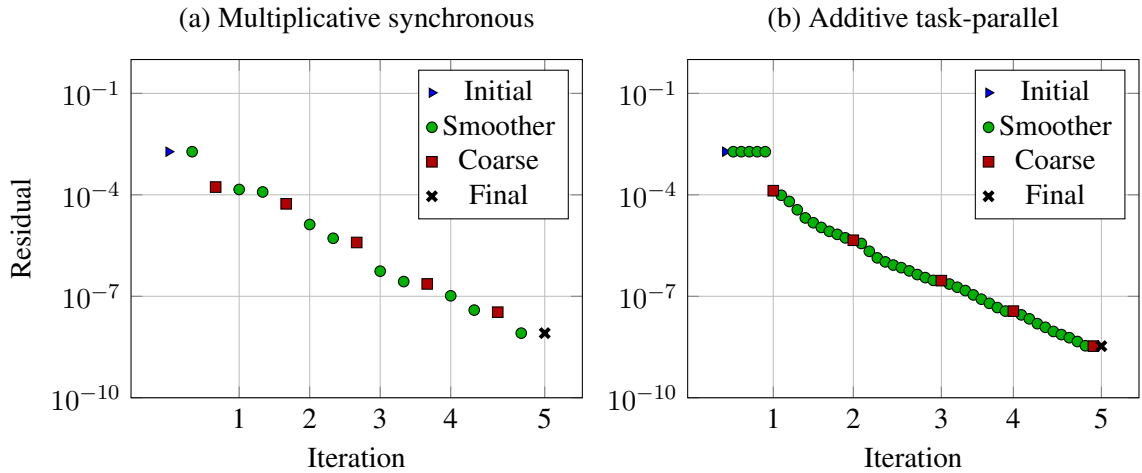
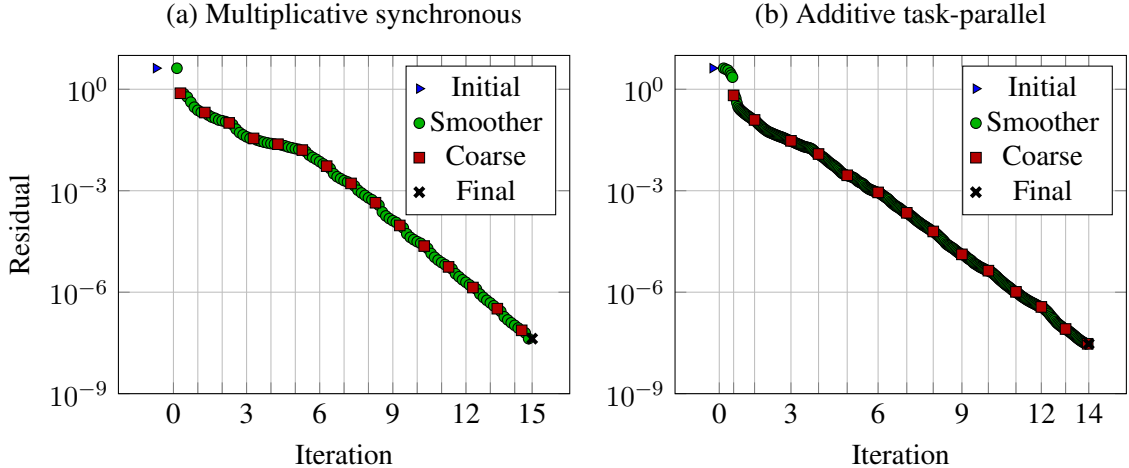
\begin{figure}[ht!]
    \centering
    \begin{tikzpicture}
\begin{groupplot}[
  group style={
    group size=2 by 1,
    horizontal sep=14mm, vertical sep=16mm,
    xlabels at=edge bottom,
    ylabels at=edge left,
  },
  ymode=log,
  xlabel={Iteration},
  ylabel={Residual},
  grid=both,
  legend pos=north east,
  ymin=1e-9, ymax=1e+1,      
  scale only axis, width=6cm, height=4.5cm,
]


\pgfkeys{
  /mods/.is family,
  /mods/initialMod/.initial=1,
  /mods/coarseMod/.initial=1,
  /mods/smootherMod/.initial=1
}

\pgfkeys{
  /mods/initialMod=1,
  /mods/coarseMod=1,
  /mods/smootherMod=1
}

\newcount\initialcount
\newcount\coarsecount
\newcount\smoothercount

\initialcount=0
\coarsecount=0
\smoothercount=0

\nextgroupplot[
    title={(a) Multiplicative synchronous},
    xtick={0,7,14, 21, 28, 35, 42, 49, 56, 63, 70, 77, 84, 91, 98, 103},
    xticklabels={0,, , 3, , , 6, , , 9, , ,12 , , , 15},
]

\initialcount=0
\coarsecount=0
\smoothercount=0

\addplot+[
  only marks,
  mark=triangle*,
  mark size=2pt,
  mark options={fill=blue!70!blue,draw=black,rotate=-90},
  x filter/.code={%
    \IfStrEq{\thisrow{type}}{initial}{
        \advance\initialcount by 1
        \ifnum\initialcount=\pgfkeysvalueof{/mods/initialMod}
            \initialcount=0
        \else
            \def\pgfmathresult{}
        \fi
    }{
        \def\pgfmathresult{}
    }
  },
] table[x=step,y=residual,col sep=comma]{figures/stokesCore48_Res128_p3_Mutliplicative.txt};
\addlegendentry{Initial}

\addplot+[
  only marks,
  mark=*,
  mark size=2pt,
  mark options={fill=green!70!black,draw=black},
  x filter/.code={%
    \IfStrEq{\thisrow{type}}{smoother}{
        \advance\smoothercount by 1
        \ifnum\smoothercount=\pgfkeysvalueof{/mods/smootherMod}
            \smoothercount=0
        \else
            \def\pgfmathresult{}
        \fi
    }{
        \def\pgfmathresult{}
    }
  },
] table[x=step,y=residual,col sep=comma]{figures/stokesCore48_Res128_p3_Mutliplicative.txt};
\addlegendentry{Smoother}

\addplot+[
  only marks,
  mark=square*,
  mark size=2pt,
  mark options={fill=red!70!black,draw=black},
  x filter/.code={%
    \IfStrEq{\thisrow{type}}{coarse}{
        \advance\coarsecount by 1
        \ifnum\coarsecount=\pgfkeysvalueof{/mods/coarseMod}
            \coarsecount=0
        \else
            \def\pgfmathresult{}
        \fi
    }{
        \def\pgfmathresult{}
    }
  },
] table[x=step,y=residual,col sep=comma]{figures/stokesCore48_Res128_p3_Mutliplicative.txt};
\addlegendentry{Coarse}

\addplot+[
  only marks,
  mark=x,
  mark size=2.5pt,
  mark options={draw=black!70!black, line width=1.5pt},
  x filter/.code={%
    \IfStrEq{\thisrow{type}}{final}{}{%
      \def\pgfmathresult{}
    }%
  },
] table[x=step,y=residual,col sep=comma]{figures/stokesCore48_Res128_p3_Mutliplicative.txt};
\addlegendentry{Final}


\nextgroupplot[title={(b) Additive task-parallel},
    xtick={-6,7,26, 59, 84, 110, 138, 164, 188, 214, 238, 263, 286, 308, 329},
    xticklabels={0,, , 3, , , 6, , , 9, , ,12 , , 14},
]

\initialcount=0
\coarsecount=0
\smoothercount=0

\addplot+[
  only marks,
  mark=triangle*,
  mark size=2pt,
  mark options={fill=blue!70!blue,draw=black,rotate=-90},
  x filter/.code={%
    \IfStrEq{\thisrow{type}}{initial}{
        \advance\initialcount by 1
        \ifnum\initialcount=\pgfkeysvalueof{/mods/initialMod}
            \initialcount=0
        \else
            \def\pgfmathresult{}
        \fi
    }{
        \def\pgfmathresult{}
    }
  },
] table[x=step,y=residual,col sep=comma]{figures/stokesCore48_Res128_p3_TaskParallel.txt};
\addlegendentry{Initial}

\addplot+[
  only marks,
  mark=*,
  mark size=2pt,
  mark options={fill=green!70!black,draw=black},
  x filter/.code={%
    \IfStrEq{\thisrow{type}}{smoother}{
        \advance\smoothercount by 1
        \ifnum\smoothercount=\pgfkeysvalueof{/mods/smootherMod}
            \smoothercount=0
        \else
            \def\pgfmathresult{}
        \fi
    }{
        \def\pgfmathresult{}
    }
  },
] table[x=step,y=residual,col sep=comma]{figures/stokesCore48_Res128_p3_TaskParallel.txt};
\addlegendentry{Smoother}

\addplot+[
  only marks,
  mark=square*,
  mark size=2pt,
  mark options={fill=red!70!black,draw=black},
  x filter/.code={%
    \IfStrEq{\thisrow{type}}{coarse}{
        \advance\coarsecount by 1
        \ifnum\coarsecount=\pgfkeysvalueof{/mods/coarseMod}
            \coarsecount=0
        \else
            \def\pgfmathresult{}
        \fi
    }{
        \def\pgfmathresult{}
    }
  },
] table[x=step,y=residual,col sep=comma]{figures/stokesCore48_Res128_p3_TaskParallel.txt};
\addlegendentry{Coarse}

\addplot+[
  only marks,
  mark=x,
  mark size=2.5pt,
  mark options={draw=black!70!black, line width=1.5pt},
  x filter/.code={%
    \IfStrEq{\thisrow{type}}{final}{}{%
      \def\pgfmathresult{}
    }%
  },
] table[x=step,y=residual,col sep=comma]{figures/stokesCore48_Res128_p3_TaskParallel.txt};
\addlegendentry{Final}

\end{groupplot}
\end{tikzpicture}
    \caption{Convergence comparison of configurations of OrthoMG for the 2D Stokes problem with the additive Schwarz smoother, polynomial degree $p=2$, resolution of $128^2$ cells simulated with 384 cores \label{fig:convCompStokes}}
\end{figure}

This section compares the convergence behavior of the proposed task-parallel OrthoMG method to the multiplicative synchronous variant.
We consider the same simulations as in \cref{sec:smootherComp}, focusing on the additive Schwarz smoother as it demonstrated superior behavior.

Figure~\ref{fig:convComp} and Figure~\ref{fig:convCompStokes} illustrate the residual decay over iterations for the \emph{multiplicative synchronous} and \emph{additive task-parallel} variants of the selected example cases for benchmark problems. 
Both variants display a monotone reduction in the residual norm due to the orthonormalization and minimization process inherent to the OrthoMG method.
The largest single-step decreases occur during the initial coarse-grid corrections, which rapidly reduce the residual norm by several orders of magnitude.
This behavior can be attributed to the well-known advantage of multigrid methods with effective elimination of low-frequency error components. 
Subsequent smoothing steps dominate the intermediate iterations, maintaining a near-geometric decay over steps with a regular trend.

The task-parallel variant shows a brief initial lag due to the overlapping execution,
 as the coarse-grid correction is performed with the stale data until the first synchronization.
However, it quickly compensates for this lag as it performs more smoothing steps within each outer cycle,
 even leading to a faster overall convergence in terms of wall-clock time with approximately 2.4 seconds in comparison to 4.1 seconds (cf., \cref{fig:smootherComp}) for the Poisson problem.
For the Stokes problem, both configurations exhibit nearly identical convergence rates and runtimes (5.4 vs 5.6 seconds), with the task-parallel approach achieving convergence in slightly fewer iterations (14 vs. 15).
Here, the difference in the number of smoother steps with iterations becomes more nuanced due to the large number of cores, 
 indicating that the task-parallel approach effectively balances the trade-off between asynchronicity and convergence efficiency.
Due to the runtime adjustment discussed in \cref{sec:adaptiveLoops} in the task-parallel configuration,
 its number of smoothing steps per cycle is not fixed but varies dynamically.
This feature provides a significant advantage to the task-parallel approach, as it allows more smoothing sweeps when the coarse-grid correction is delayed.

Tables \ref{tab:poisson2d_schwarz_blockjacobi_384} and \ref{tab:stokes2d_schwarz_blockjacobi_384}
 present the iteration counts for two solver configurations with varying polynomial degrees and grid resolutions.
In all tested configurations, the additive task-parallel variant consistently exhibit a comparative number of iterations to reach convergence compared to the multiplicative synchronous scheme.
This indicates that the task-parallel execution does not compromise convergence efficiency despite its lagged semi-asynchronous nature.
Moreover, in some cases, the task-parallel approach even outperforms the classical variant in terms of iteration counts,
 which can be attributed to the increased number of smoothing steps performed per cycle due to the adaptive loops.
 This is result of the overlapping execution, which allows the smoother to proceed independently and perform additional iterations while waiting for coarse-grid corrections on dedicated cores,
 which in comparison can be idle or inefficient in the classical synchronous scheme.
Although theoretically, the semi-asynchronous execution could lead to less effective error reduction per iteration due to stale data usage,
 in practice, the results suggest that the increased smoothing can compensate for this potential drawback.
Moreover, the overall runtime benefits of the task-parallel approach, as demonstrated in \cref{sec:runtimeComp}, further showcase its effectiveness.

\subsection{Scalability of different variants} \label{sec:runtimeComp}
\begin{figure}[ht!]
    \centering
    \pgfplotscreateplotcyclelist{solverstyles}{
  {blue, mark=o,        mark options={scale=1, fill=white, line width=0.9pt}, thick}, 
  {blue, mark=*,        mark options={scale=1},                                    thick}, 
  {red,  mark=triangle, mark options={scale=1.5, fill=white, line width=0.9pt}, thick}, 
  {red,  mark=triangle*,mark options={scale=1.5},                                    thick}, 
  {green!50!black, mark=square,  mark options={scale=1, fill=white, line width=0.9pt}, thick}, 
  {green!50!black, mark=square*, mark options={scale=1},                                    thick}, 
  {black, dashed, ultra thick}, 
}

\pgfplotsset{
  myStrongScaling/.style={
    xlabel={Cores ($3 \times 2^k$)},
    ylabel={Wall-clock time [s]},
    xmode=log, ymode=log, log basis x=2,
    xmin=24, xmax=6144,
    ymin=1,  ymax=3200,
    xtick={48,96,192,384,768, 1536, 3072},
    xticklabels={$2^{4}$,$2^{5}$,$2^{6}$,$2^{7}$,$2^{8}$,$2^{9}$,$2^{10}$},
    xmajorgrids, ymajorgrids, xminorgrids, yminorgrids,
    major grid style={opacity=0.3},
    minor grid style={opacity=0.15},
    cycle list name=solverstyles,
    legend cell align=left,
    scale only axis, width=6cm, height=4.5cm,
    after end axis/.code={
      \node[anchor=north east, inner sep=0pt, xshift=12pt, yshift=-3.5pt]
        at (rel axis cs:0,0) {$3\times$};
    },
  }
}

\newlength{\legsampleshift}
\setlength{\legsampleshift}{2.2em} 

\pgfplotsset{
  legend image post style={xshift=\legsampleshift}, 
  legend style={/tikz/column sep=0pt, nodes={inner xsep=0pt}}, 
}

\newcommand{\legendpad}[1]{\makebox[\legsampleshift][l]{#1}} 

\begin{tikzpicture}

\begin{groupplot}[
  group style={group size=2 by 2, horizontal sep=14mm, vertical sep=16mm,  xlabels at=edge bottom,  ylabels at=edge left},
  myStrongScaling
]

\nextgroupplot[title={(a) Additive synchronous}, legend to name=SharedLegend, legend columns=2, legend style={
  font=\footnotesize,
  /tikz/every even column/.append style={column sep=0.9em},
  transpose legend=true,
} ]

\addplot coordinates
  { (48,14.0674496) (96,12.3093008) (192,20.0229718)
    (384,11.2530691) (768,10.8203978) (1536,32.2772909)
    (3072,37.8639002) };
\addlegendentry{$p=2, \, 32^3$}

\addplot coordinates {
  (48,47.17369045) (96,31.9123323) (192,37.6712404)
  (384,33.9358699) (768,42.8408877) (1536,84.49125115)
  (3072,120.6794786)
};
\addlegendentry{$p=3,\,32^3$}

\addplot coordinates
  { (48,208.4390613) (96,110.6308991)
    (192,64.92741865) (384,57.88264955)
    (768,41.39015035) (1536,106.2201375)
    (3072,68.120944) };
\addlegendentry{$p=2, \, 48^3$}

\addplot coordinates {
  (48,651.2375134) (96,288.5578546) (192,162.6761327)
  (384,168.1187838) (768,96.80848325)
  (1536,228.1174124) (3072,215.2006908)
};
\addlegendentry{$p=3,\,48^3$}

\addplot coordinates
  { (48,131.7574022) (96,169.459177)
    (192,126.514921) (384,58.4898594)
    (768,50.280057) (1536,276.4195734)
    (3072,417.5012465) };
\addlegendentry{$p=2, \, 64^3$}

\addplot coordinates {
  (192,345.145479) (384,441.1811891)
  (768,1253.184774)  (1536,384.5817082)
  (3072,1060.6818) 
};
\addlegendentry{$p=3,\,64^3$}

\pgfmathsetmacro{\PrefA}{96}
\pgfmathsetmacro{\TrefA}{150}
\addplot[dashed, thick, samples at={48,96,192,384,768, 1536, 3072}] ({x},{\TrefA*\PrefA/x});
\addlegendentry{Ideal scaling}

\nextgroupplot[title={(b) Multiplicative synchronous},]

\addplot coordinates
  { (48,4.7583463) (96,4.4509757)
    (192,6.1322959) (384,3.6276339)
    (768,3.75580035) (1536,8.4999716)
    (3072,12.45108385) };

\addplot coordinates
  { (48,12.1471932) (96,21.5619375)
    (192,6.1446) (384,6.81701305)
    (768,9.235057) (1536,13.3018942)
    (3072,21.93692445) };

\addplot coordinates
  { (48,31.07597985) (96,27.8978166)
    (192,19.78534825) (384,15.08857615)
    (768,12.8631392) (1536,19.8322156)
    (3072,23.02138215) };

\addplot coordinates
  { (48,65.520957) (96,58.40848735)
    (192,47.93871665) (384,71.8909073)
    (768,35.6575659) (1536,38.13251985)
    (3072,22.4070564) };

\addplot coordinates
  { (48,48.82269565) (96,33.4844666)
    (192,30.120824) (384,16.2398956)
    (768,16.2427642) (1536,33.78376675)
    (3072,95.17654505) };

\addplot coordinates
  { (96,64.0879037) (192,45.69461705)
    (384,41.876316) (768,60.41838105)
    (1536,28.40703745) (3072,120.3343898) };

\pgfmathsetmacro{\PrefB}{96}\pgfmathsetmacro{\TrefB}{150}
\addplot[dashed, thick, samples at={48,96,192,384,768, 1536, 3072}, forget plot] ({x},{\TrefB*\PrefB/x});

\nextgroupplot[title={(c) Additive task-parallel}]

\addplot coordinates
  { (48,9.30420835) (96,7.80204515)
    (192,4.76241375) (384,3.23084005)
    (768,2.35796355) (1536,2.56184515)
    (3072,3.09432705) };

\addplot coordinates
  { (48,15.33634955) (96,11.06741175)
    (192,6.6449238) (384,4.84297575)
    (768,3.5047155) (1536,3.09555525)
    (3072,3.6326855) };

\addplot coordinates
  { (48,45.19439085) (96,31.7290998)
    (192,16.9766386) (384,13.30731785)
    (768,7.3824886) (1536,4.70402515)
    (3072,4.02765015) };

\addplot coordinates
  { (96,60.9931935) (192,26.5797173)
    (384,19.0524432) (768,10.8903418)
    (1536,8.40954395) (3072,5.72071775) };

\addplot coordinates
  { (96,92.11040315) (192,44.79867005)
    (384,23.0110218) (768,14.92874645)
    (1536,10.3197651) (3072,5.6801967) };

\addplot coordinates
  { (384,35.70371295) (768,20.7797149)
    (1536,12.83091175) (3072,11.48938145) };

\pgfmathsetmacro{\PrefC}{96}\pgfmathsetmacro{\TrefC}{150}
\addplot[dashed, thick, samples at={48,96,192,384,768, 1536, 3072}, forget plot] ({x},{\TrefC*\PrefC/x});

\nextgroupplot[title={(d) Hybrid (Multiplicative-Task-Parallel)}, ylabel={}]

\addplot coordinates
  { (48,7.49580575) (96,5.1652778)
    (192,5.00387335) (384,3.97502835)
    (768,3.45228715) (1536,4.94813175)
    (3072,9.1484312) };

\addplot coordinates
  { (48,11.5842654) (96,9.6639572)
    (192,6.2169642) (384,6.69978415)
    (768,5.0943657) (1536,10.512147)
    (3072,10.07086395) };

\addplot coordinates
  { (48,36.70964515) (96,24.0420807)
    (192,13.6039344) (384,15.76503105)
    (768,12.62044365) (1536,8.5374343)
    (3072,15.060647) };

\addplot coordinates
  { (48,68.15882465) (96,48.87939165)
    (192,28.13786695) (384,18.1924186)
    (768,13.39583575) (1536,12.05397055)
    (3072,14.77885855) };

\addplot coordinates
  { (48,41.7996265) (96,32.64128605)
    (192,22.0420909) (384,13.38510385)
    (768,11.8765239) (1536,12.4073917)
    (3072,10.87145985) };

\addplot coordinates
  { (96,47.90919275) (192,35.71111665)
    (384,20.0891015) (768,11.6176562)
    (1536,11.78255845) (3072,12.9143584) };

\pgfmathsetmacro{\PrefD}{96}\pgfmathsetmacro{\TrefD}{150}
\addplot[dashed, thick, samples at={48,96,192,384,768, 1536, 3072}, forget plot] ({x},{\TrefD*\PrefD/x});

\end{groupplot}

\end{tikzpicture}
\begin{center}
\hspace{1cm} \pgfplotslegendfromname{SharedLegend}
\end{center}
    \caption{Runtime comparison of configurations of OrthoMG for the 3D Poisson problem with varying number of cores, polynomial degree and resolution \label{fig:runTimeCompPoisson}}
\end{figure}
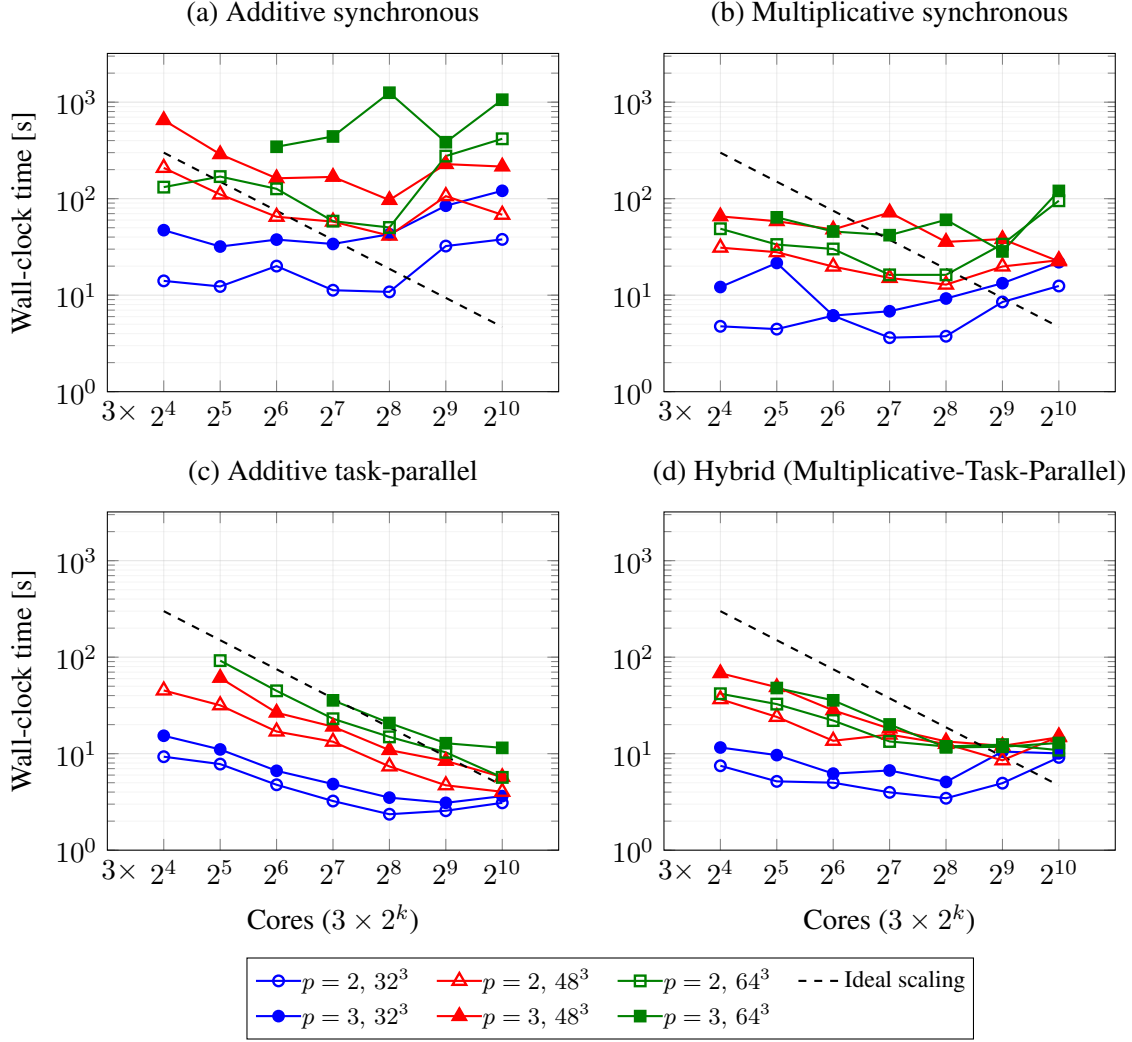

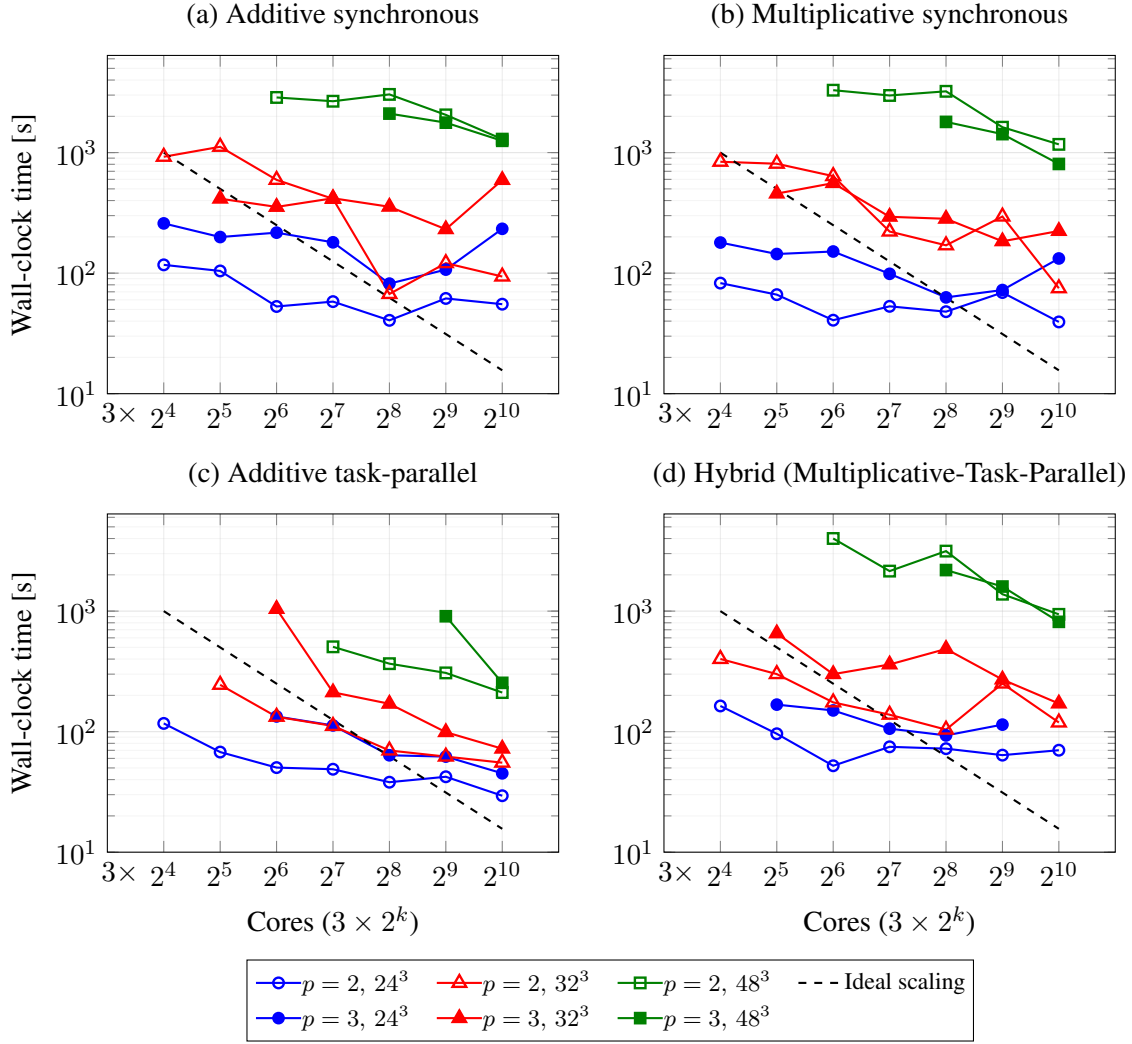
\begin{figure}[ht!]
    \centering
    \pgfplotscreateplotcyclelist{solverstyles}{
  {blue, mark=o,        mark options={scale=1, fill=white, line width=0.9pt}, thick}, 
  {blue, mark=*,        mark options={scale=1},                                    thick}, 
  {red,  mark=triangle, mark options={scale=1.5, fill=white, line width=0.9pt}, thick}, 
  {red,  mark=triangle*,mark options={scale=1.5},                                    thick}, 
  {green!50!black, mark=square,  mark options={scale=1, fill=white, line width=0.9pt}, thick}, 
  {green!50!black, mark=square*, mark options={scale=1},                                    thick}, 
  {black, dashed, ultra thick}, 
}

\pgfplotsset{
  myStrongScaling/.style={
    xlabel={Cores ($3 \times 2^k$)},
    ylabel={Wall-clock time [s]},
    xmode=log, ymode=log, log basis x=2,
    xmin=24, xmax=6144,
    ymin=10,  ymax=6400,
    xtick={48,96,192,384,768, 1536, 3072},
    xticklabels={$2^{4}$,$2^{5}$,$2^{6}$,$2^{7}$,$2^{8}$,$2^{9}$,$2^{10}$},
    xmajorgrids, ymajorgrids, xminorgrids, yminorgrids,
    major grid style={opacity=0.3},
    minor grid style={opacity=0.15},
    cycle list name=solverstyles,
    legend cell align=left,
    scale only axis, width=6cm, height=4.5cm,
    after end axis/.code={
      \node[anchor=north east, inner sep=0pt, xshift=12pt, yshift=-3.5pt]
        at (rel axis cs:0,0) {$3\times$};
    },
  }
}

\setlength{\legsampleshift}{2.2em} 

\pgfplotsset{
  legend image post style={xshift=\legsampleshift}, 
  legend style={/tikz/column sep=0pt, nodes={inner xsep=0pt}}, 
}

\newcommand{\legendpad}[1]{\makebox[\legsampleshift][l]{#1}} 

\begin{tikzpicture}
\pgfmathsetmacro{\PrefA}{48}
\pgfmathsetmacro{\TrefA}{1000}

\begin{groupplot}[
  group style={group size=2 by 2, horizontal sep=14mm, vertical sep=16mm,  xlabels at=edge bottom,  ylabels at=edge left},
  myStrongScaling
]

\nextgroupplot[title={(a) Additive synchronous}, legend to name=SharedLegend, legend columns=2, legend style={
  font=\footnotesize,
  /tikz/every even column/.append style={column sep=0.9em},
  transpose legend=true,
} ]

\addplot coordinates
  { (48,117.1691992) (96,104.2988202) (192,52.89512945)
    (384,58.0582189) (768,40.597543)
    (1536,61.67971395) (3072,55.26798905) };

\addlegendentry{$p=2, \, 24^3$}

\addplot coordinates {
  (48,258.6653685) (96,199.3043191) (192,216.7110396)
  (384,180.2270318) (768,81.8690132)
  (1536,106.9792946) (3072,233.2009538)
};
\addlegendentry{$p=3,\, 24^3$}

\addplot coordinates
  { (48,923.4865398) (96,1117.571359) (192,594.2040018)
    (384,413.4739137) (768,67.16526293)
    (1536,121.0377029) (3072,93.7141556) };
\addlegendentry{$p=2, \, 32^3$}

\addplot coordinates {
  (96,414.590164) (192,354.3344402) (384,418.0846136)
  (768,355.3622497) (1536,231.3648166) (3072,592.2773854)
};
\addlegendentry{$p=3,\, 32^3$}

\addplot coordinates
  { (192,2867.574152) (384,2667.795033)
    (768,3039.860235) (1536,2055.853519) (3072,1293.23683) };
\addlegendentry{$p=2, \, 48^3$}

\addplot coordinates {
  (768,2106.33589)
  (1536,1772.500992) (3072,1252.406857)
};
\addlegendentry{$p=3, \, 48^3$}

\addplot[dashed, thick, samples at={48,96,192,384,768, 1536, 3072}] ({x},{\TrefA*\PrefA/x});
\addlegendentry{Ideal scaling}

\nextgroupplot[title={(b) Multiplicative synchronous},]

\addplot coordinates
  { (48,82.70713195) (96,66.3222307) (192,40.700433)
    (384,53.05222425) (768,47.90567465)
    (1536,69.08420985) (3072,39.3860403) };

\addplot coordinates {
  (48,179.4452751) (96,144.1135605) (192,151.2208383)
  (384,98.74106335) (768,62.92678045)
  (1536,72.43230965) (3072,132.3574315)
};

\addplot coordinates
  { (48,839.6685301) (96,808.6642553) (192,638.2228713)
    (384,221.4161306) (768,170.2093237)
    (1536,294.8817014) (3072,74.95365135) };

\addplot coordinates {
  (96,457.0432054) (192,558.459041) (384,293.3899642)
  (768,282.4778733) (1536,183.7570922) (3072,223.7090325)
};

\addplot coordinates
  { (192,3294.491394) (384,2978.5794) (768,3224.409262)
    (1536,1624.427768) (3072,1170.614036) };

\addplot coordinates {
  (768,1798.117795) (1536,1423.461236) (3072,805.1137144)
};

\addplot[dashed, thick, samples at={48,96,192,384,768, 1536, 3072}] ({x},{\TrefA*\PrefA/x});

\nextgroupplot[title={(c) Additive task-parallel},]

\addplot coordinates
  { (48,117.0225981) (96,67.7206644) (192,50.32092605)
    (384,48.77096325) (768,38.07229145)
    (1536,42.23478435) (3072,29.4039833) };

\addplot coordinates {
  (192,133.6536679) (384,112.6115232) (768,63.664166)
  (1536,62.08887645) (3072,45.1021811)
};

\addplot coordinates
  { (96,244.7369351) (192,132.8818335) (384,111.5289655)
    (768,69.7860221) (1536,62.10502908)
    (3072,55.29200835)  };

\addplot coordinates {
  (192,1041.643672) (384,211.9952041) (768,170.2862691)
  (1536,99.02943448) (3072,72.19678395)
};

\addplot coordinates
  { (384,505.5956017) (768,366.4745105)
    (1536,307.1546812) (3072,211.0818766) };

\addplot coordinates {
  (1536,906.4059114) (3072,254.056479)
};

\addplot[dashed, thick, samples at={48,96,192,384,768, 1536, 3072}] ({x},{\TrefA*\PrefA/x});

\nextgroupplot[title={(d) Hybrid (Multiplicative-Task-Parallel)}, ylabel={}]

\addplot coordinates
  { (48,163.3504713) (96,96.0255562) (192,52.11536665)
    (384,74.8548189) (768,72.24478905)
    (1536,63.84063575) (3072,70.1348196) };

\addplot coordinates {
  (96,167.5908226) (192,149.9718225) (384,105.878163)
  (768,93.05016645) (1536,114.3616786)
};

\addplot coordinates
  { (48,400.8348611) (96,300.2597268) (192,174.9362273)
    (384,138.1842305) (768,103.9588123)
    (1536,250.7264202) (3072,119.1122496) };

\addplot coordinates {
  (96,653.5964054) (192,299.8602103) (384,360.8602933)
  (768,485.4280322) (1536,270.6092038) (3072,171.0528719)
};

\addplot coordinates
  { (192,4003.998923) (384,2147.580706) 
    (768,3146.593121) (1536,1381.01905) (3072,939.9981396) };

\addplot coordinates {
   (768,2190.476064)
  (1536,1599.090869) (3072,814.5347301)
};

\addplot[dashed, thick, samples at={48,96,192,384,768, 1536, 3072}] ({x},{\TrefA*\PrefA/x});

\end{groupplot}

\end{tikzpicture}
\begin{center}
\hspace{1cm} \pgfplotslegendfromname{SharedLegend}
\end{center}
    \caption{Runtime comparison of configurations of OrthoMG for the 3D Stokes problem with varying number of cores, polynomial degree and resolution \label{fig:runTimeCompStokes}}
\end{figure}

To assess the efficiency of the different solver configurations, their runtimes are compared for 3D Poisson and Stokes problems with discontinuous coefficients, as described in \cref{sec:benchmarks}. 
These test cases are identical to the 2D configuration discussed in previous sections, differing only in dimensionality.
As the smoother, the additive Schwarz method is employed, since it demonstrated superior performance in \cref{sec:smootherComp}.

For the Poisson problem, the computational grid is resolved with $32\times32\times32$, $48\times48\times48$ and $64\times64\times64$ cells.
The polynomial degree is varied between $p=2$ and $p=3$, resulting in $10$ or $20$ DOFs per cell, which corresponds to an exact doubling of the total number of DOFs.
The total number of DOFs ranges from approximately $0.34$ million to $5.3$ million, including the extended cut cells.

The Stokes problem is resolved with $24\times24\times24$, $32\times32\times32$ and $48\times48\times48$ cells.
The polynomial degree is varied between $p=2$ and $p=3$, resulting in $34$ or $70$ DOFs per cell.
The corresponding total number of DOFs ranges from approximately $0.48$ million to $7.8$ million.

\cref{fig:runTimeCompPoisson} and \cref{fig:runTimeCompStokes} compare the wall-clock time on a log-log scale with an ideal $T\!\propto\!1/$cores reference for both problems. 
As can be seen, all configurations display reduced times with increased core counts initially, but the degree of scalability varies significantly between them.
Unsurprisingly, the worst scalability behavior is observed with \emph{the additive synchronous} configuration,
 which suffers from excessive synchronization overheads and communication bottlenecks at higher core counts with intrinsically poorer convergence.
This leads even to runtime increases in some cases, indicating a breakdown of scalability.
A similar situation is also observed with \emph{the multiplicative synchronous} strategy, which, despite being initially competitive, 
degrades rapidly at higher core counts due to coarse-level bottlenecks and synchronization delays.

On the other hand, the hybrid scheme exhibits a competitive regime over the advantageous task-parallel configuration for the relatively low number of cores ($\leq 384$) for the Poisson problem, 
 but its scalability degrades at higher core counts.
This indicates that even a single level of classical multiplicative execution can hinder scalability.

Across all problem sizes and polynomial degrees, 
\emph{the additive task-parallel} configuration exhibits the most favorable strong-scaling, for example,
for $64^3$ with $p{=}2$, the time decreases from $92.1$\,s at 96 cores to $14.9$\,s at 768 cores, 
corresponding to a speed-up of $\sim\!6.2$ over $8\times$ cores, with an efficiency of approximately $77\%$
but then flattens at the highest concurrencies due to insufficient workload per core to $5.7$\,s at 3072 cores,
 corresponding to a speed-up of $\sim\!16.2$ over $32\times$ cores.
It should also be noted that, the maximum memory constraint of $25\times10^3$ DOFs per core limits the problem size for the given core counts, 
 while smaller problems cannot be executed at the highest core counts efficiently due to insufficient workload per core.

Another observation related to the synchronicity is the sensitivity to problem decomposition and load balancing. 
As can be seen on the presented figures, the runtime does not always monotonously decrease for the synchronous schemes, even before reaching the saturation points,
 while the semi-asynchronous task-parallel variant exhibits relatively consistent reductions. 
Although somehow speculative, this can be attributed to the stronger sensitivity of synchronous schemes to local load imbalances and communication hotspots, 
 where any overloaded subdomain or unfavorable partitioning forces all ranks to wait, leading to irregular runtime behavior.
In contrast, the task-parallel approach can mitigate the impact of these imbalances and result in a more regular behavior.
This is particularly important for problems with strong anisotropies or heterogeneous coefficients as in the presented multiphase benchmarks.

As presented, the task-parallel configuration consistently outperforms the other variants for sufficiently large problems and core counts,
 demonstrating superior scalability and the best runtimes. 
However, it also comes with certain disadvantages as memory requirement per core increases due to the assignment of the fine level grids to fewer cores. 
This ultimately limits the minimum problem size that can be addressed on a given hardware setup, in comparison to the other configurations.
Moreover, Figures~\ref{fig:runTimeCompPoisson} and \ref{fig:runTimeCompStokes} exhibit that with relatively small problem sizes or low number of cores, 
the performance differences between the configurations are less pronounced or can even be worse for the task-parallel scheme due to its additive, overlapped scheme.

Overall, the results indicate that minimizing global synchronization to a single point per cycle and
 concurrent smoother and coarse-grid solves are essential to translate additional cores into reduced runtime. 
When these conditions are met (e.g., by the task-parallel variant), 
the scaling behavior remains close to the ideal reference over a wide core range and only flattens at the highest concurrencies, 
where the problem size becomes too small for the allocated number of cores. 
In contrast, configurations with either strictly sequential phases or synchronous overlap exhibit premature saturation or even runtime regressions,
 highlighting the sensitivity of multigrid solvers to communication/computation balance on modern nodes.
In addition, it is noted that having a runtime adjustment between smoother and coarse-grid groups helps to improve the convergence behavior,
 particularly when combined with the additive Schwarz smoother, which shows robust performance across various problem sizes and polynomial degrees. 

\section{Conclusion and outlook} \label{sec:Conclusion}
We have presented a semi-asynchronous task-parallel orthonormalization multigrid method for solving large-scale linear systems
 arising from the discretization of partial differential equations.
The proposed method employs task-based parallelism to temporally overlap smoothing and coarse-grid correction operations,
 reducing idle times and improving resource utilization on distributed-memory systems.
The presented numerical results demonstrate that the task-parallel configuration outperforms traditional synchronous approaches in terms of runtime and strong scalability,
    particularly at high numbers of cores.
This is due to the effective reduction of synchronization bottlenecks and coarse-level delays, leading to significant speed-ups on high core counts.

Moreover, the presented method utilizes a runtime adjustment between smoother and coarse-grid groups,
 allowing for adaptive balancing of their workloads and even improved convergence rates.
It is observed that this is particularly useful when combined with the additive Schwarz smoother,
 which exhibits robust performance across various problem sizes and polynomial degrees.
However, the method also introduces limitations, such as increased memory requirement due to the assignment of grids to fewer cores.
This may restrict the minimum problem size that can be solved on a given hardware setup compared to other configurations.
In addition, the overlapped execution may lead to initial convergence lags due to stale data usage until the first synchronization occurs.

As future work, these limitations could be addressed by exploring more memory-efficient data structures and adaptive synchronization strategies to mitigate the impact of stale data.
It could also be of interest to investigate further optimizations, such as asynchronous smoothers or convergence criteria based on problem characteristics.
Furthermore, investigating the effect of load balancing and coarsening strategies on the performance of the task-parallel method could reveal further improvements.
In addition, extending the method to more complex PDEs and heterogeneous computing environments with GPUs could broaden its applicability and enhance its performance in practical scenarios.

\section*{Acknowledgments}
The work of T. Toprak is funded by the Deutsche Forschungsgemeinschaft (DFG, German Research Foundation) – Project-ID 492661287 – TRR 361. 
Additionally, the work of the authors is supported by the Graduate
School CE within the Centre for Computational Engineering at
TU Darmstadt. 
Furthermore, the authors gratefully acknowledge the computing time provided to them on the high-performance computer Lichtenberg II at TU Darmstadt, 
funded by the German Federal Ministry of Education and Research (BMBF) and the State of Hesse.

\section*{Data availability}
The data that support the findings of this study are available at DOI and can be reproduced by the presented worksheets with the open-source BoSSS framework.
\section*{Conflict of interest}
The authors declare no potential conflict of interests.


\bibliographystyle{unsrtnat}
\bibliography{template}  






\end{document}